\documentclass[letterpaper, paper,11pt]{AAS}		

\usepackage{bm}
\usepackage{amsmath}
\usepackage{amssymb}
\usepackage{subfigure}
\usepackage[colorlinks=true, pdfstartview=FitV, linkcolor=black, citecolor= black, urlcolor= black]{hyperref}
\usepackage{overcite}
\usepackage{footnpag}			      	
\usepackage{booktabs}
\usepackage[hang,flushmargin]{footmisc}

\PaperNumber{21-651}

\begin{document}

\title{Computation and Analysis of Jupiter-Europa and Jupiter-Ganymede Resonant Orbits in the Planar Concentric Circular Restricted 4-Body Problem}

\author{Bhanu Kumar\thanks{PhD Candidate, School of Mathematics, Georgia Institute of Technology, 686 Cherry St. NW, Atlanta, GA 30332.},  
Rodney L. Anderson\thanks{Technologist, Jet Propulsion Laboratory, California Institute of Technology, 4800 Oak Grove Drive, Pasadena, CA 91109.},
Rafael de la Llave\thanks{Professor, School of Mathematics, Georgia Institute of Technology, 686 Cherry St. NW, Atlanta, GA 30332.},
\ and Brian Gunter\thanks{Associate Professor, Daniel Guggenheim School of Aerospace Engineering, Georgia Institute of Technology, 270 Ferst Drive, Atlanta, GA 30332, AIAA Associate Fellow}
}

\maketitle{}

\begin{abstract}
Many unstable periodic orbits of the planar circular restricted 3-body problem (PCRTBP) persist as invariant tori when a periodic forcing is added to the equations of motion. In this study, we compute tori corresponding to exterior Jupiter-Europa and interior Jupiter-Ganymede PCRTBP resonant periodic orbits in a concentric circular restricted 4-body problem (CCR4BP). Motivated by the 2:1 Laplace resonance between Europa and Ganymede's orbits, we then attempt the continuation of a Jupiter-Europa 3:4 resonant orbit from the CCR4BP into the Jupiter-Ganymede PCRTBP. We strongly believe that the resulting dynamical object is a KAM torus lying near but not on the 3:2 Jupiter-Ganymede resonance. \end{abstract}

\section{Introduction}

\footnotetext{\copyright \, 2021. All rights reserved.}

Many studies using dynamical systems theory and incorporating the use of resonant orbits for tour design have focused on the tour endgame, or the final few resonances before approaching a moon. These studies generally use as their model the circular restricted three-body problem (CRTBP), which incorporates the effects of a primary and secondary body. One particular problem that has been of interest has been the final approach to Europa after a series of flybys of the Galilean moons using either ballistic, impulsive, or low-thrust trajectories\cite{Anderson2010, Anderson2011, Anderson2011b}. Our recent work has enabled the rapid and accurate computation of resonant orbits and heteroclinic connections to quickly design trajectories that traverse these resonances within both the CRTBP and the elliptic RTBP \cite{kumar2021journal, kumar2020, kumar2021feb}.

One challenge that has recently arisen is the design of trajectories that transfer from flybys of Ganymede to flybys of Europa\cite{Anderson2021b}. These flybys correspond to resonant periodic orbits in the Jupiter-Ganymede and Jupiter-Europa CRTBP models. This problem has been approached in the past by using various approximations, such as patched-conic models\cite{Sweetser1997Trajectory} or patching two CRTBP models together\cite{Anderson2021a}. To more accurately model these types of trajectories, a full four-body model may be utilized. Earlier work has explored some planar Lyapunov periodic orbits within a concentric model\cite{blazevski2012}.

In our study, we focus on the computation of unstable resonant orbits known to be useful for the final transition between Ganymede and Jupiter in a Jupiter-Ganymede-Europa-Spacecraft model. Specifically, we use the planar concentric circular restricted four-body problem (CCR4BP) as the dynamical model for this analysis. This model assumes that the moons travel in concentric circular orbits around Jupiter, and is an example of a periodically-forced planar CRTBP (PCRTBP) model; in fact, the CCR4BP can be considered a periodic perturbation of \emph{two different CRTBP models}, the Jupiter-Europa and Jupiter-Ganymede models. It is known\cite{kumar2021rapid} that when a periodic forcing is applied to a CRTBP, up to some forcing size, most unstable periodic orbits from that CRTBP persist as 2D unstable quasi-periodic orbits (also known as invariant tori). It is these tori that we compute. 

One property peculiar to the Jupiter-Europa-Ganymede system is that Europa and Ganymede are themselves in an approximate 2:1 mean motion resonance with each other, known as the Laplace resonance\cite{Barnes2011}; Europa makes approximately two revolutions around Jupiter in the time that Ganymede makes one. Hence, a spacecraft orbit resonant with either of Europa or Ganymede is also close to being resonant with the other. Motivated by this, after continuing CRTBP unstable Jupiter-Europa resonant periodic orbits into the Jupiter-Europa-Ganymede CCR4BP, we then also try to continue the resulting CCR4BP quasi-periodic orbits into the Jupiter-Ganymede CRTBP.

\section{Planar Concentric Circular Restricted 4-body Problem} \label{modelsection}
The planar concentric circular restricted 4-body problem (CCR4BP) models the motion of a spacecraft under the gravitational influence of three large bodies of masses $m_{1}$, $m_{2}$, and $m_{3}$. It is assumed that $m_{1} >> m_{2}, m_{3}$, as is the case when $m_{1}$ is Jupiter and $m_{2}$ and $m_{3}$ are Europa and Ganymede. In the concentric model\cite{blazevski2012} considered here, $m_{2}$ and $m_{3}$ are assumed to revolve in concentric circles around $m_{1}$ of radii $r_{12}$ and $r_{13}$, with no effect of $m_{2}$ on the motion of $m_{3}$ and vice versa. The angular velocities $\Omega_{2}$ and $\Omega_{3}$ of the revolution of $m_{2}$ and $m_{3}$ around $m_{1}$ can be found using Kepler's third law, and are given by
\begin{equation} \label{kepler} \Omega_{i} = \sqrt{\frac{\mathcal{G}(m_{1}+m_{i})}{r_{1i}^{3}}} , \quad i = 2,3\end{equation}
where $\mathcal{G}$ is the universal gravitational constant. This is not a coherent model; the motion of $m_{1}$, $m_{2}$, and $m_{3}$ just described is not a solution of the full 3-body problem with finite masses, but it is a very good approximation of their true motion. We consider the planar CCR4BP, so it is assumed that both circular orbits as well as the motion of the spacecraft all lie in the same plane. 

Now, define mass ratios $\mu = \frac{m_{2}}{m_{1}+ m_{2}}$ and  $\mu_{3} = \frac{m_{3}}{m_{1}+ m_{2}}$. Next, normalize length, mass, and time units so that $r_{12}$ becomes 1, $\mathcal{G}(m_{1}+m_{2})$ becomes $1$, and thus $\Omega_{2}$ also becomes $1$. We can then write the planar CCR4BP equations of motion for the spacecraft in a synodic coordinate system exactly similar to the usual rotating coordinate frame for the CRTBP. In particular, $m_{1}$ and $m_{2}$ lie on the $x$-axis of this synodic frame, and the origin is at the $m_{1}$-$m_{2}$ barycenter. In this frame and units, the angle between the position of $m_{3}$ and the $x$-axis at time $t$ is $\theta_{3}(t) = (\Omega_{3}-1)t + \theta_{3,0}$, so the position of $m_{3}$ is $(x_{3}(t), y_{3}(t)) = (-\mu+r_{13} \cos(\theta_{3}), r_{13} \sin(\theta_{3}) )$. The equations of motion are then given in position-momentum space $(x,y,p_{x},p_{y})$ by (see Blazevski and Ocampo\cite{blazevski2012} for a derivation)
\begin{equation} \label{ccr4bpEOM} \begin{gathered} \dot x= p_{x} + y \quad \quad \dot y= p_{y} - x \\
 \dot p_x=  p_{y} - (1-\mu)\frac{x+\mu}{r_{1}^{3}} - \mu \frac{x - (1-\mu)}{r_{2}^{3}} - \mu_{3} \frac{x - x_{3}}{r_{3}^{3}} - \mu_{3} \frac{\cos \theta_{3}}{r_{13}^{2}} \\
 \dot p_y=  -p_{x} - (1-\mu)\frac{y}{r_{1}^{3}} - \mu \frac{y}{r_{2}^{3}} - \mu_{3} \frac{y - y_{3}}{r_{3}^{3}} - \mu_{3} \frac{\sin \theta_{3}}{r_{13}^{2}} \end{gathered} \end{equation}
where $r_{1} = \sqrt{(x+\mu)^{2} + y^{2}}$, $r_{2} = \sqrt{(x-1+\mu)^{2} + y^{2}}$, and $r_{3} = \sqrt{(x-x_{3})^{2} + (y-y_{3})^{2}}$ are the distances from the spacecraft to $m_{1}$, $m_{2}$, and $m_{3}$, respectively.  Note that we recover the usual $m_{1}$-$m_{2}$ planar CRTBP equations of motion if we set $\mu_{3}=0$. The equations of motion given by Equation \eqref{ccr4bpEOM} are Hamiltonian, with time-periodic Hamiltonian function
\begin{equation}  \label{ccr4bpH} H_{\mu_{3}}(x,y,p_x,p_{y}, \theta_{3})= \frac{p_{x}^{2}+p_{y}^{2}}{2} + p_{x}y -p_{y}x - \frac{1-\mu}{r_{1}} - \frac{\mu}{r_{2}}  - \frac{\mu_{3}}{r_{3}} + \mu_{3} \frac{x \cos \theta_{3}}{r_{13}^{2}} + \mu_{3} \frac{y \sin \theta_{3}}{r_{13}^{2}} \end{equation}

Just as we normalized units to make $r_{12}=\mathcal{G}(m_{1}+m_{2})=\Omega_{2}=1$ and wrote the planar CCR4BP equations of motion in an $m_{1}$-$m_{2}$ synodic frame, we can alternatively normalize units to make $r_{13}=\mathcal{G}(m_{1}+m_{3})=\Omega_{3}=1$, and then write the equations of motion in an $m_{1}$-$m_{3}$ synodic frame centered at the $m_{1}$-$m_{3}$ barycenter. In this case, the equations of motion are of the same form as Equation \eqref{ccr4bpEOM}; defining $\bar \mu = \frac{m_{3}}{m_{1}+ m_{3}}$ and  $\bar \mu_{2} = \frac{m_{2}}{m_{1}+ m_{3}}$, the $m_{1}$-$m_{3}$ rotating frame equations are
\begin{equation} \label{ccr4bpEOM_frame3} \begin{gathered} \dot x= p_{x} + y \quad \quad \dot y= p_{y} - x \\
 \dot p_x=  p_{y} - (1-\bar \mu)\frac{x+\bar \mu}{r_{1}^{3}} - \bar \mu \frac{x - (1-\bar \mu)}{r_{3}^{3}} - \bar \mu_{2} \frac{x - x_{2}}{r_{2}^{3}} - \bar \mu_{2} \frac{\cos \theta_{2}}{r_{12}^{2}} \\
 \dot p_y=  -p_{x} - (1-\bar \mu)\frac{y}{r_{1}^{3}} - \bar \mu \frac{y}{r_{3}^{3}} - \bar \mu_{2} \frac{y - y_{2}}{r_{2}^{3}} - \bar \mu_{2} \frac{\sin \theta_{2}}{r_{12}^{2}} \end{gathered} \end{equation}
where $r_{1} = \sqrt{(x+\bar \mu)^{2} + y^{2}}$, $r_{2} = \sqrt{(x-x_{2})^{2} + (y-y_{2})^{2}}$, and $r_{3} = \sqrt{(x-1+ \bar \mu)^{2} + y^{2}}$ are the distances from the spacecraft to $m_{1}$, $m_{2}$, and $m_{3}$, respectively, and $\theta_{2}(t) = (\Omega_{2}-1)t + \theta_{2,0}$, $(x_{2}(t), y_{2}(t)) = (-\bar \mu + r_{12} \cos(\theta_{2}), r_{12} \sin(\theta_{2}) )$.  We recover the $m_{1}$-$m_{3}$ planar CRTBP equations of motion when $\bar \mu_{2}=0$. As both Eq. \eqref{ccr4bpEOM} and \eqref{ccr4bpEOM_frame3} model the same system, there exists a transformation which takes a state from the $m_{1}$-$m_{2}$ frame to the $m_{1}$-$m_{3}$ frame. We present this next. 

\subsection{Transformation from $m_{1}$-$m_{2}$ frame to $m_{1}$-$m_{3}$ frame} \label{transformsection}
We know that the $m_{1}$-$m_{3}$ frame $x$-axis is the line between those two bodies, which is represented in the $m_{1}$-$m_{2}$ coordinate frame as a line passing through the point $(-\mu, 0)$, making an angle $\theta_{3}$ with the $m_{1}$-$m_{2}$ frame $x$-axis, and revolving at an angular rate of $\Omega_{3}-1$. Thus, to transform a state from the $m_{1}$-$m_{2}$ frame into the $m_{1}$-$m_{3}$ frame described earlier, one must follow several steps:

\begin{enumerate} 
\item Shift the origin of the position from the $m_{1}$-$m_{2}$ barycenter to $m_{1}$
\item Rotate the $m_{1}$-$m_{2}$ frame position and velocity vectors by the angle $-\theta_{3}$
\item Apply the transport theorem to get the actual $m_1$-$m_{3}$ frame apparent velocity 
\item Rescale the length, mass, and time units, and set $\theta_{2} = -\theta_{3}$
\item Shift the origin of the position from $m_{1}$ to the $m_{1}$-$m_{3}$ barycenter
\end{enumerate}

Let $(x,y, p_{x}, p_y, \theta_{3})$ be a state in the $m_{1}$-$m_{2}$ frame, and let $(\bar x, \bar y,  \bar p_x, \bar p_y, \theta_{2})$ denote the equivalent state in the $m_{1}$-$m_{3}$ frame. Denote $R_{\rho}$ as the $2 \times 2$ rotation matrix such that for $\bold{x} \in \mathbb{R}^{2}$, the vector $R_{\rho} \bold{x}$ is rotated by $\rho$ radians counterclockwise. Then, remembering that $\dot x = p_{x}+y$ and $\dot y=p_{y}-x$, steps 1 and 2 from above for finding $(\bar x, \bar y, \bar p_{x}, \bar p_{y})$ require computing 
\begin{equation} \begin{bmatrix} x'  \\ y' \end{bmatrix}
 = R_{-\theta_{3}} \begin{bmatrix} x+\mu \\ y \end{bmatrix} \quad \quad  \begin{bmatrix} \dot x'  \\ \dot y' \end{bmatrix}
 = R_{-\theta_{3}} \begin{bmatrix} \dot x \\ \dot y \end{bmatrix}  \end{equation}
 After this, we can write expressions for $(\bar x, \bar y, \dot{\bar  x}, \dot{\bar  y})$ in the $m_{1}$-$m_{3}$ frame as
\begin{equation} \label{m1m3_state_vels}  
\bar x = \frac{x'}{r_{13}} - \bar \mu \quad \quad \bar y = \frac{y'}{r_{13}} \quad \quad \dot{\bar x} = \frac{\dot x' +(\Omega_{3}-1)y'}{r_{13} \Omega_{3}}  \quad \quad \dot{\bar y} = \frac{\dot y' -(\Omega_{3}-1)x'}{r_{13} \Omega_{3}} 
  \end{equation}
These expressions encompass steps 3, 4, and 5 from the previous list. Last of all, we can compute $\bar p_x = \dot{\bar x}- \bar y$ and $\bar p_{y}=\dot{\bar y} + \bar x$; substituting Eq. \eqref{m1m3_state_vels} into these expressions for $\bar p_x$ and $\bar p_y$ actually gives an even simpler alternative to Eq. \eqref{m1m3_state_vels}. In particular, we find that
\begin{equation} \label{m1m3_state}  
\bar p_x = \frac{\dot x' -y'}{r_{13} \Omega_{3}}  = \frac{p_x'}{r_{13} \Omega_{3}} \quad \quad \quad \bar p_y = \frac{\dot y' + x'}{r_{13} \Omega_{3}}  - \bar \mu = \frac{p_y'}{r_{13} \Omega_{3}}  - \bar \mu 
  \end{equation}
  where $\begin{bmatrix} p_x'  & p_y' \end{bmatrix}^{T} = R_{-\theta_{3}} \begin{bmatrix} p_x & p_y+\mu \end{bmatrix}^{T} $. Finally, once  $(\bar x, \bar y,  \bar p_x, \bar p_y)$ has been found, one should also keep in mind that $\theta_{2} = - \theta_{3}$, which completes the full extended state required to propagate the $m_{1}$-$m_{3}$ equations of motion given in Eq. \eqref{ccr4bpEOM_frame3}. We will refer to this transformation as $\Phi_{\theta_{3}}:\mathbb{R}^{4}\rightarrow \mathbb{R}^{4}$; that is, using the earlier notation of this section, $(\bar x, \bar y,  \bar p_x, \bar p_y)=\Phi_{\theta_{3}}(x,y, p_{x}, p_y)$. Note that $\Phi_{\theta_{3}}$ is just a composition of affine and linear transformations, so its derivative $D\Phi_{\theta_{3}}$ is simple to compute. 
 
\section{Computing Invariant Tori in the Planar CCR4BP} \label{backgroundsection}

The planar CCR4BP is an example of a periodically-perturbed PCRTBP model. As mentioned earlier, for $\mu_{3} = 0$ in Eq. \eqref{ccr4bpEOM} and $\bar \mu_{2} = 0$ in Eq. \eqref{ccr4bpEOM_frame3}, the CCR4BP reduces to either the $m_{1}$-$m_{2}$ or the $m_{1}$-$m_{3}$ CRTBP, respectively. On the other hand, for $\mu_{3} > 0$ ($\bar \mu_{2}>0$), Eq. \eqref{ccr4bpEOM} (Eq, \eqref{ccr4bpEOM_frame3}) consists of the planar $m_{1}$-$m_{2}$ ($m_{1}$-$m_{3}$) CRTBP equations of motion plus perturbation terms dependent on the angle $\theta_{3}$ ($\theta_{2}$), which advances at the constant rate $\dot \theta_{3} = \Omega_{3}-1$ ($\dot \theta_{2} = \Omega_{2}-1$). The perturbation terms are periodic, with period $T_{p} = \frac{2\pi}{|\Omega_3-1|}$ ($\bar T_{p} = \frac{2\pi}{|\Omega_2-1|}$). Furthermore, as evidenced by Eq. \eqref{ccr4bpH} (a similar Hamiltonian exists for Eq. \eqref{ccr4bpEOM_frame3}), the additional perturbative terms due to $\mu_{3}$ or $\bar \mu_{2}$ do not break the Hamiltonian structure of the equations of motion. 

Due to this structure of the planar CCR4BP as a Hamiltonian periodic-perturbation of both the $m_{1}$-$m_{2}$ planar CRTBP as well the $m_{1}$-$m_{3}$ planar CRTBP, some conclusions can be drawn. First of all, it can be concluded that up to some value of the perturbation parameter $\mu_{3}$ or $\bar \mu_{2}$, unstable periodic orbits from both $m_{1}$-$m_{2}$ and $m_{1}$-$m_{3}$ CRTBP's will generally persist as 2D invariant tori in the CCR4BP extended phase space $(x,y,p_{x},p_{y}, \theta_{3})$ or $(\bar x,\bar y,\bar p_{x},\bar p_{y}, \theta_{2})$. This persistence is related to that of normally hyperbolic invariant manifolds and of Kolmogorov-Arnold-Moser (KAM) tori, as is explained in our previous paper \cite{kumar2021rapid}. For any resulting 2D torus, one frequency will be that of the perturbation ($\dot \theta_{3}$ or $\dot \theta_{2}$), while the other will be the frequency of the original CRTBP periodic orbit. As we will show, an invariant 2D torus in the $m_{1}$-$m_{2}$ synodic CCR4BP coordinates can be transformed to one in the $m_{1}$-$m_{3}$ rotating frame using the transformation $\Phi_{\theta_{3}}$. 

In earlier work\cite{kumar2021rapid}, we developed efficient methods for computing unstable invariant tori and their center, stable, and unstable directions in periodically-perturbed planar CRTBP models. To apply these methods to the planar CCR4BP, we first need the concept of stroboscopic maps. 

\subsection{Stroboscopic Maps} \label{stroboscopic}

The quasi-periodic orbits of interest in the planar CCR4BP lie on 2D unstable invariant tori in the 5D extended phase space $(x,y,p_{x},p_{y},\theta_{p})$, where $p=3$ or $2$ depending on the frame being used. These invariant tori can be parameterized as the image of a function of two angles $K_2 : \mathbb{T}^{2} \rightarrow \mathbb{R}^{4} \times \mathbb{T}$. Any quasi-periodic trajectory $\bold{x}(t)$ lying on this torus can be expressed as
	\begin{equation} \label{torusXTraj} \bold{x}(t) = K_2(\theta, \theta_p)  \quad \quad \quad \theta = \theta_{0}+\Omega_1 t, \quad  \theta_p = \theta_{p,0}+(\Omega_p-1) t \end{equation}
where $\bold{x}(0)$ determines $\theta_{0}$ and $\theta_{p,0}$. $\theta_{p}$ and $\Omega_p-1$ are the perturbation phase angle and frequency described earlier, respectively. Defining the stroboscopic map $F: \mathbb{R}^{4} \times \mathbb{T} \rightarrow \mathbb{R}^{4} \times \mathbb{T}$ as the time-$\frac{2\pi}{|\Omega_p-1|}$ mapping of extended phase space points by the equations of motion, we have
	 \begin{equation} \label{invariance_fK2} F(K_2(\theta, \theta_p)) = K_2(\theta+\omega, \theta_p), \text{ where  } \omega = 2\pi \Omega_1/|\Omega_p-1| \end{equation} 
since the angle $\theta_p$ increases or decreases by $2\pi$ in the time $\frac{2\pi}{|\Omega_p-1|}$. Since the value of $\theta_{p}$ is invariant under the map $F$, we can fix $\theta_p$ and define $K(\theta) = K_2(\theta, \theta_p)$. Without loss of generality, we choose $\theta_{p}=0$ in this study; this has the added benefit that both $\theta_{2}=\theta_{3}=0$, so the $m_{1}$-$m_{2}$ and $m_{1}$-$m_{3}$ frames are aligned with each other at this phase. With $K$ thus defined, Equation \eqref{invariance_fK2} becomes 
\begin{equation} \label{invariance} F(K(\theta)) = K(\theta+\omega)  \end{equation} 
Ignoring the constant $\theta_{p}$ component of the extended phase space and making a slight abuse of notation, we have $F:\mathbb{R}^{4} \rightarrow \mathbb{R}^{4}$ and $K:\mathbb{T} \rightarrow \mathbb{R}^{4}$. Equation \eqref{invariance} implies that $K$ is an invariant 1D torus (circle) of $F$. Hence, basing our computations on the stroboscopic map $F$ is more efficient than solving for tori invariant under the flow of the equations of motion, since we reduce the phase space dimension from 5D to 4D and the dimension of the unknown invariant tori from 2D to 1D. Once the 1D stroboscopic map torus is computed, one can visualize the full 2D torus of the flow by numerically integrating points from the invariant circle by the equations of motion Eq. \eqref{ccr4bpEOM} or \eqref{ccr4bpEOM_frame3}.

\subsection{Parameterization Methods for Tori and Bundles} \label{quasiNewton}

With the stroboscopic map $F$ defined by the CCR4BP equations of motion, we now wish to find solutions $K(\theta)$ of Equation \eqref{invariance}. In this section, we summarize the results of the parameterization methods developed in Kumar, Anderson, and de la Llave\cite{kumar2021rapid} for computing invariant tori in periodically-perturbed planar CRTBP models such as the planar CCR4BP; these methods were themselves  inspired by methods described in Haro et al. \cite{haroetal} and Zhang and de la Llave\cite{zhang}. The rotation number $\omega=\frac{2\pi \Omega_1}{|\Omega_p-1|}$ is generally known; for instance, if the CCR4BP torus being solved for comes from a known PCRTBP periodic orbit, then $\Omega_{1}$ is the frequency of that periodic orbit. 

In our previous work\cite{kumar2021rapid} , we developed an efficient  quasi-Newton method for solving equation \eqref{invariance} given a sufficiently accurate initial guess. Our quasi-Newton method adds an extra equation to be solved in addition to Eq. \eqref{invariance}.  In particular, as well as $K(\theta)$, we simultaneously solve for matrix-valued periodic functions $P(\theta)$, $\Lambda(\theta): \mathbb{T} \rightarrow \mathbb{R}^{4 \times 4}$ satisfying
		 \begin{equation}  \label{bundleEquations} DF(K(\theta)) P(\theta) = P(\theta+\omega) \Lambda(\theta) \end{equation} 
$P(\theta)$ and $\Lambda(\theta)$ are the matrices of bundles and of Floquet stability, respectively. For each $\theta \in \mathbb{T}$, the columns of $P(\theta)$ are comprised of the (linearly independent) tangent, symplectic conjugate center, stable, and unstable directions of the torus at the point $K(\theta)$, in that order, while $\Lambda$ has the form
\begin{equation}  \label{Lambdaform}
\Lambda(\theta)=\begin{bmatrix}
1 &  T(\theta)   & 0 & 0 \\ 0 &  1   & 0 & 0 \\ 0 & 0  & \lambda_s  & 0 \\ 0 &  0 & 0 & \lambda_u \end{bmatrix}
 \end{equation}
where $T:\mathbb{T} \rightarrow \mathbb{R}$ and $\lambda_s, \lambda_u \in \mathbb{R}$ are constants with $\lambda_s<1$ and $\lambda_u > 1$. $\lambda_s$ and $\lambda_u$ are the stable and unstable multipliers for the invariant torus, which we will also find useful in this study. 

As it turns out, solving simultaneously for $K$, $P$, and $\Lambda$  not only gives stability information, but actually has lower computational complexity than more commonly-used methods which solve for $K$ alone. When discretizing functions of $\theta$ at $N$ evenly spaced angle values, our method requires only $O(N)$ storage and $O(N \log N)$ operations since the nearly-diagonal $\Lambda$ matrix allows decoupling of the scalar functional equations which are solved in each quasi-Newton step\cite{kumar2021rapid}. This is much more efficient than $K$-only methods which require solving an un-decoupled system of functional equations in each step, which after discretization requires $O(N^{3})$ operations due to Gaussian elimination.

\subsection{Continuation of Tori into the Planar CCR4BP} \label{quasiNewton}

In this work, we use the quasi-Newton method of the previous section as part of a numerical continuation scheme for computing tori. Note that the stroboscopic map $F$ defined earlier is obtained by integrating points by the equations of motion Eq. \eqref{ccr4bpEOM} or \eqref{ccr4bpEOM_frame3}; hence, $F$ depends on the value of $\mu_{3}$ or $\bar \mu_{2}$, depending on the frame being used. Due to this, in this section we write $F=F_{\mu_{3}}$ or $F=F_{\bar \mu_{2}}$ to signify this dependence on parameters which we will vary during the continuation. 

To compute an invariant torus for some desired value of $\mu_{3}$, we can start with a periodic orbit in the $m_{1}$-$m_{2}$ PCRTBP, and then continue this by $\mu_{3}$ until the required $\mu_{3}$ value is reached. A periodic orbit is an invariant circle for $F_{\mu_{3}=0}$, so we can take this periodic orbit as a starting point for the continuation by $\mu_{3}$. We can similarly continue $m_{1}$-$m_{3}$ PCRTBP periodic orbits by $\bar \mu_{2}$. One thing to note is that due to Kepler's third law (see Eq. \eqref{kepler}), fixing the angular velocity $\Omega_{3}$ of $m_{3}$ and varying $\mu_{3}$ during the continuation results in $r_{13}$ changing as well; similarly, $r_{12}$ changes during the continuation by $\bar \mu_{2}$ in the case of an $m_{1}$-$m_{3}$ frame continuation. It is desirable to fix $\Omega_{3}$ or $\Omega_{2}$ to its physical value during the entire continuation, so that the invariant circle rotation number $\omega=\frac{2\pi \Omega_1}{|\Omega_3-1|}$ or $\frac{2\pi \Omega_1}{|\Omega_2-1|}$ remains constant. Once a single torus has been continued by $\mu_{3}$ or $\bar \mu_{2}$ to the desired physical mass value, the same quasi-Newton method can also be used to continue the resulting torus by $\omega$ (with fixed mass values) to find other CCR4BP tori in the same family. 

While the quasi-Newton method gives a way to correct an approximate solution of Eq. \eqref{invariance}-\eqref{bundleEquations}, we improved continuation performance by adding a predictor step based on the Poincar\'e-Lindstedt method to complement the corrector. We described the Lindstedt method for $\omega$ continuation in our previous work\cite{kumar2021rapid}. Here we derive a similar method for predicting the next torus during a continuation by $\mu_{3}$ (the predictor for $\bar \mu_{2}$ continuation works the exact same way). 

Assume we have a solution $K_{\mu_{3}}(\theta), P_{\mu_{3}}(\theta), \Lambda_{\mu_{3}}(\theta)$ to Eq.\eqref{invariance}-\eqref{bundleEquations} for some value of $\mu_{3}$. Then, we have that $F_{\mu_{3}}(K_{\mu_{3}}(\theta)) = K_{\mu_{3}}(\theta+\omega)$. Differentiating this with respect to $\mu_{3}$, we have
\begin{equation} \label{lindstedt} \frac{dF_{\mu_{3}}}{d\mu_{3}}(K_{\mu_{3}}(\theta)) + DF_{\mu_{3}}(K_{\mu_{3}}(\theta)) \frac{dK_{\mu_{3}}}{d\mu_{3}}(\theta)= \frac{dK_{\mu_{3}}}{d\mu_{3}}(\theta+\omega) \end{equation}
where $DF_{\mu_{3}}$ is the matrix-valued derivative of $F_{\mu_{3}}$ with respect to the state $(x,y,p_{x},p_{y})$, and $\frac{dF_{\mu_{3}}}{d\mu_{3}}$ is the vector-valued derivative of $F_{\mu_{3}}$ with respect to the parameter $\mu_{3}$. To solve for $\frac{dK_{\mu_{3}}}{d\mu_{3}}(\theta)$, substitute $\frac{dK_{\mu_{3}}}{d\mu_{3}}(\theta) = P_{\mu_{3}}(\theta) \xi(\theta)$ into Eq. \eqref{lindstedt}. Left multiplying by $P^{-1}_{\mu_{3}}(\theta+\omega)$ and recalling from Eq. \eqref{bundleEquations} that $P_{\mu_{3}}^{-1}(\theta+\omega)DF_{\mu_{3}}(K_{\mu_{3}}(\theta)) P_{\mu_{3}}(\theta) =  \Lambda_{\mu_{3}}(\theta) $, we get the equation
\begin{equation} \label{lindstedt2}  \Lambda_{\mu_{3}}(\theta) \xi(\theta) - \xi(\theta+\omega) = -P_{\mu_{3}}^{-1}(\theta+\omega)\frac{dF_{\mu_{3}}}{d\mu_{3}}(K_{\mu_{3}}(\theta))  \end{equation}
Recalling that $F_{\mu_{3}}$ just represents the integration of an ODE by some fixed time, $\frac{dF_{\mu_{3}}}{d\mu_{3}}(K_{\mu_{3}}(\theta))$ can be computed using the usual variational equations for a parameter dependent system. Hence, the RHS of Eq. \eqref{lindstedt2} is known; equations of the form of Eq. \eqref{lindstedt2} are almost diagonal due to Eq. \eqref{Lambdaform}, and $\xi$ can easily be solved for using Fourier series as is described in our previous work\cite{kumar2021rapid}. Once $\xi$ and hence $\frac{dK_{\mu_{3}}}{d\mu_{3}}(\theta)$ are known, we can use this to predict the torus for $\mu_{3} + \Delta \mu_{3}$ as
\begin{equation} \label{predictor} K_{\mu_{3}+ \Delta \mu_{3} }(\theta)  \approx K_{\mu_{3}}(\theta) + \Delta \mu_{3}\frac{dK_{\mu_{3}}}{d\mu_{3}}(\theta) \end{equation}
The result of this computation is then used as the initial guess (along with $P_{\mu_{3}}(\theta)$ and $\Lambda_{\mu_{3}}(\theta)$) for the quasi-Newton method for finding an invariant circle of $F_{\mu_{3} + \Delta \mu_{3}}$. 

\subsection{Transforming Tori and Bundles Between Frames} \label{quasiNewton}

Assume we have a solution $K, P, \Lambda$ of Eq. \eqref{invariance}-\eqref{bundleEquations} with $F=F_{\mu_{3}}$, the stroboscopic map for some fixed $\theta_{3}$ of the $m_{1}$-$m_{2}$ frame equations of motion Eq. \eqref{ccr4bpEOM}; we would like to transform this to the corresponding solution $\bar K, \bar P, \bar \Lambda$ of Eq. \eqref{invariance}-\eqref{bundleEquations} for $F = F_{\bar \mu_{2}}$, the $\theta_{2}=-\theta_{3}$ based stroboscopic map of the $m_{1}$-$m_{3}$ frame equations of motion Eq. \eqref{ccr4bpEOM_frame3}. Recall the function $\Phi_{\theta_{3}}$ defined earlier taking states from the $m_{1}$-$m_{2}$ frame to the $m_{1}$-$m_{3}$ frame. The relationship between $F_{\mu_{3}}$, $F_{\bar \mu_{2}}$, and $\Phi_{\theta_{3}}$ is 
\begin{equation} \label{conjugacy1} F_{\bar \mu_{2}} \circ \Phi_{\theta_{3}} = \Phi_{\theta_{3}} \circ F_{\mu_{3}} \end{equation}
Eq. \eqref{conjugacy1} simply means that the same final state is attained whether one first transforms to the $m_{1}$-$m_{3}$ frame and then applies the $m_{1}$-$m_{3}$ frame stroboscopic map $F_{\bar \mu_{2}}$, or one first applies the $m_{1}$-$m_{2}$ frame stroboscopic map $F_{\mu_{3}}$ and then only afterwards transforms to the $m_{1}$-$m_{3}$ frame. This can be justified by noting that substituting $(x,y,p_x,p_y)=\Phi_{\theta_{3}}^{-1}(\bar x,\bar y,\bar p_x,\bar p_y)$ and $\theta_{3}=-\theta_{2}$ into Eq. \eqref{ccr4bpEOM}, followed by a rescaling of time to $m_{1}$-$m_{3}$ frame time units $\Omega_3 \,dt = d \bar t$, yields Eq. \eqref{ccr4bpEOM_frame3}; the derivation of this is simple but lengthy, so we omit it for the purpose of space. Eq. \eqref{conjugacy1} is simply the discrete-time version of the equivalency of Eq. \eqref{ccr4bpEOM} and Eq. \eqref{ccr4bpEOM_frame3} through $\Phi_{\theta_{3}}$. 

Since we know that $F_{\mu_{3}}(K(\theta)) = K(\theta+\omega)$, Eq. \eqref{conjugacy1} implies that 
\begin{equation} \label{conjugacy2}  F_{\bar \mu_{2}} \left( \Phi_{\theta_{3}} (K(\theta)) \right)  = \Phi_{\theta_{3}} \left( F_{\mu_{3}}(K(\theta) ) \right) =  \Phi_{\theta_{3}} ( K(\theta+\omega) )\end{equation}
which means that $\bar K(\theta) = \Phi_{\theta_{3}} (K(\theta))$ satisfies Eq. \eqref{invariance} for $F=F_{\bar \mu_{2}}$. Similarly, we know that $DF_{\mu_{3}}(K(\theta)) P(\theta) = P(\theta+\omega) \Lambda(\theta)$; differentiating Eq. \eqref{conjugacy1} and right-multiplying by $P$ thus gives
\begin{equation} \begin{split} \label{conjugacy2}  DF_{\bar \mu_{2}} \left( \Phi_{\theta_{3}} (K(\theta)) \right) D\Phi_{\theta_{3}} (K(\theta))  P(\theta) &= D\Phi_{\theta_{3}} \left( F_{\mu_{3}}(K(\theta) ) \right) DF_{\mu_{3}}(K(\theta) ) P(\theta) \\
&= D\Phi_{\theta_{3}} \left( K(\theta+\omega)  \right) P(\theta+\omega) \Lambda(\theta) \end{split} \end{equation}
which implies that $\bar K(\theta) = \Phi_{\theta_{3}} (K(\theta))$, $\bar P(\theta) = D\Phi_{\theta_{3}} (K(\theta)) P(\theta)$, and $\bar \Lambda(\theta) = \Lambda(\theta)$ satisfies Eq. \eqref{bundleEquations} for $F=F_{\bar \mu_{2}}$. Thus, through these equations, we can transform any solution of Eq. \eqref{invariance}-\eqref{bundleEquations} for $F = F_{ \mu_{3}}$ to the corresponding solution for $F = F_{\bar \mu_{2}}$. 

\section{Resonant Tori in the Jupiter-Europa-Ganymede Planar CCR4BP} \label{resonanttoriSection}

As mentioned in the previous sections, we can expect most $m_{1}$-$m_{2}$ and $m_{1}$-$m_{3}$ planar CRTBP unstable periodic orbits to persist as 2D invariant tori in the $m_{1}$-$m_{2}$-$m_{3}$ planar CCR4BP. In the remainder of this study, we will take $m_{1}$ to be Jupiter, $m_{2}$ to be Europa, and $m_{3}$ to be Ganymede. We know that there exist families of unstable resonant periodic orbits in the Jupiter-Europa and Jupiter-Ganymede planar CRTBP models; such a resonant periodic orbit is identified using a ratio $m$:$n$, which means that in an inertial coordinate frame, the spacecraft makes approximately $m$ revolutions around Jupiter in the same time that the moon considered makes $n$ revolutions. Using the previously described Lindstedt predictor and quasi-Newton corrector methods for continuation, we compute the invariant tori corresponding to some Jupiter-Europa and Jupiter-Ganymede resonances in the Jupiter-Europa-Ganymede CCR4BP. We refer to these as Jupiter-Europa and Jupiter-Ganymede resonant CCR4BP tori. Though initial computation of tori is done in different coordinate frames (rotating with Jupiter-Europa or Jupiter-Ganymede), we can use the transformation described earlier to visualize Jupiter-Europa resonant CCR4BP tori in a Jupiter-Ganymede frame, and vice versa. 

In all the computations to follow, the invariant tori, bundles, and multipliers were solved with an accuracy of $10^{-7}$ or better in Eq. \eqref{invariance}-\eqref{bundleEquations}. We used the values of Jupiter, Europa, and Ganymede masses and orbital periods given in Table \ref{table:constants}. The $\mathcal{G}m_{i}$ values are used to compute mass ratios; the periods are converted to the normalized time units described earlier depending on the frame being used, from which other quantities (such as $\Omega_{2}$ and $\Omega_{3}$) can be computed. 

\begin{table} 
\centering
 \begin{tabular}{c c c} 
 \hline
 Body & $\mathcal{G} m_{i}$ ($\text{m}^{3}/\text{s}^{2}$) & Orbital Period (s)  \\ [0.5ex] 
 \hline\hline
 Jupiter & $1.2668653785779600\times 10^{17}$ & -  \\ 
 \hline
 Europa & $9.8869974284299492\times 10^{12}$ & $3.0689648366400000\times 10^{5}$  \\
 \hline
 Ganymede & $3.2009998067205903\times 10^{12}$ & $6.1808096312640002\times 10^{5}$  \\
 \hline
\end{tabular}
\caption{ \label{table:constants} Masses and Orbital Periods Used for Jupiter, Europa, and Ganymede}
\end{table}

\subsection{3:4 Jupiter-Europa Resonant CCR4BP Tori}

\begin{figure}
\begin{centering}
\includegraphics[width=0.48\columnwidth]{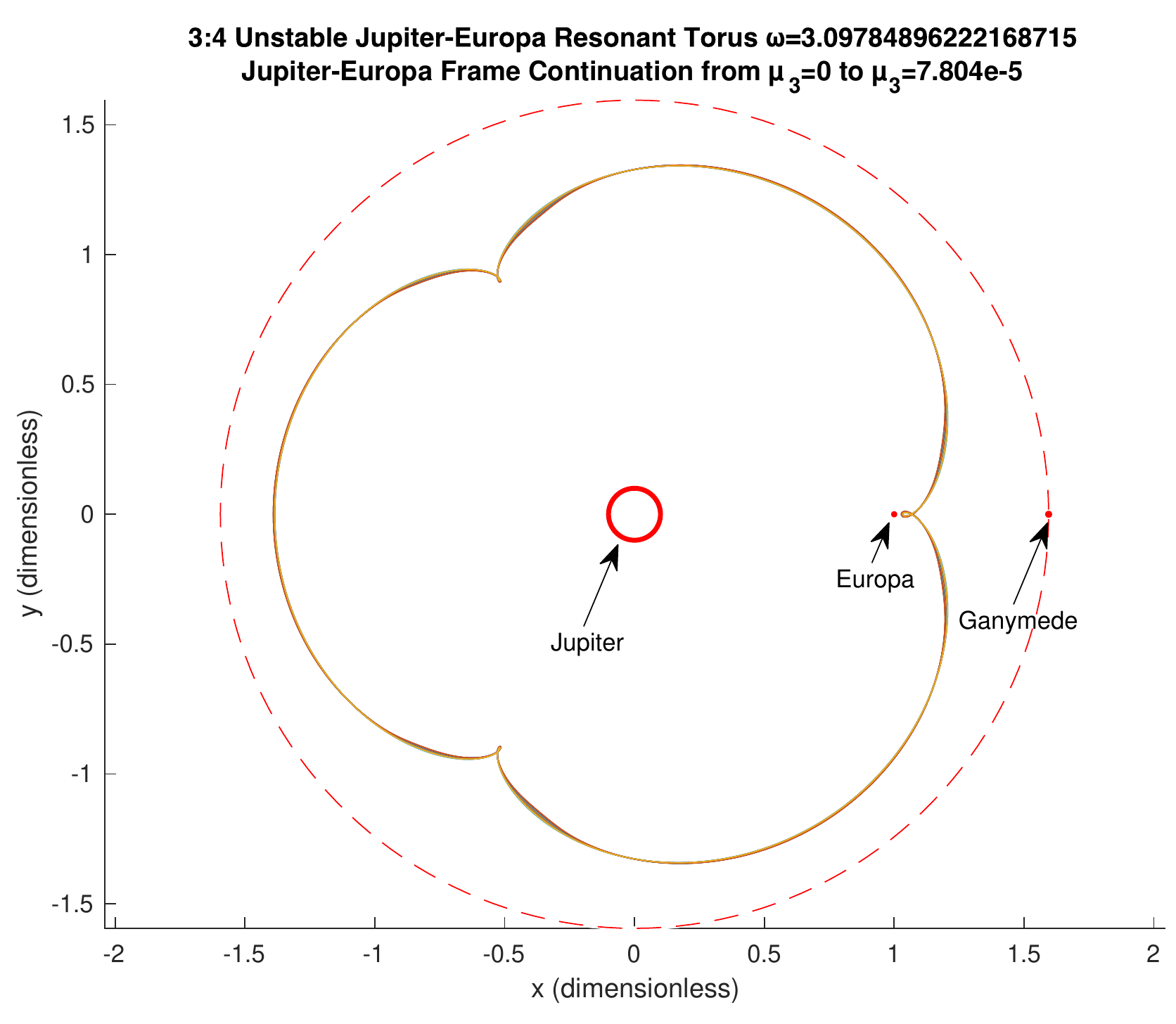}
\includegraphics[width=0.48\columnwidth]{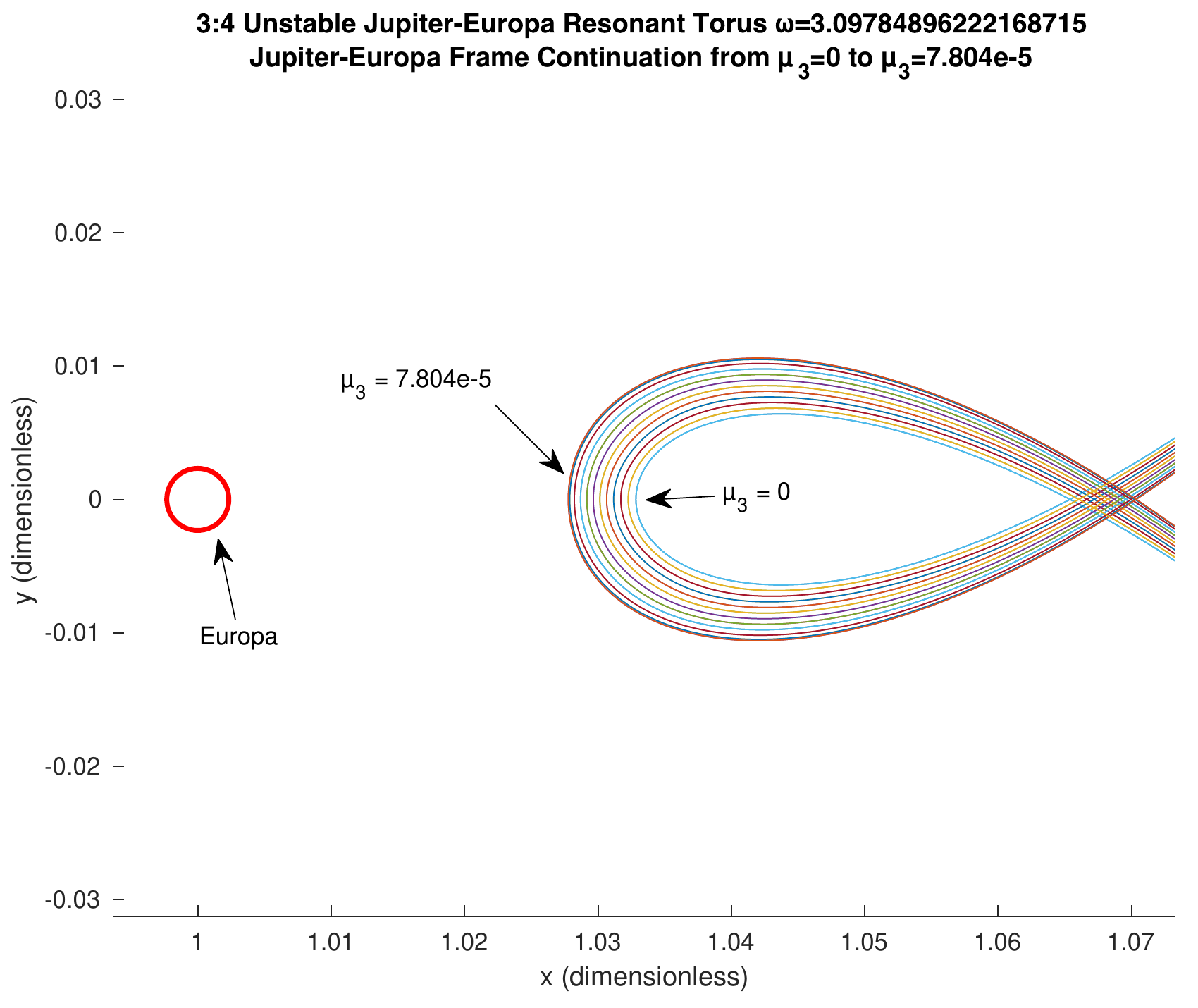}
\caption{ \label{fig:34eps} Continuation of Jupiter-Europa 3:4 Resonant CCR4BP Torus By $\mu_{3}$}
\end{centering}
\end{figure}

We first compute the quasi-periodic equivalents of Jupiter-Europa 3:4 unstable resonant periodic orbits in the CCR4BP. Starting from a periodic orbit with Jacobi constant value 3.0041 and $\omega=3.097849$ (arbitrarily chosen), the Jupiter-Europa frame planar CCR4BP stroboscopic map invariant circles computed at various $\mu_{3}$ continuation steps are displayed in Figure \ref{fig:34eps}; the continuation step size used was $\Delta \mu_{3} = 8 \times 10^{-6}$ until the final value of $\mu_{3}=7.804102777055038 \times 10^{-5}$ for Ganymede. We also show the orbit of Ganymede on the plot as a red dashed-line circle, although recall that our stroboscopic map is defined using a phase of $\theta_3=0$, so Ganymede is on the $x$-axis aligned with Europa at the phase of these computed invariant circles. 

\begin{figure}
\begin{centering}
\includegraphics[width=0.49\columnwidth]{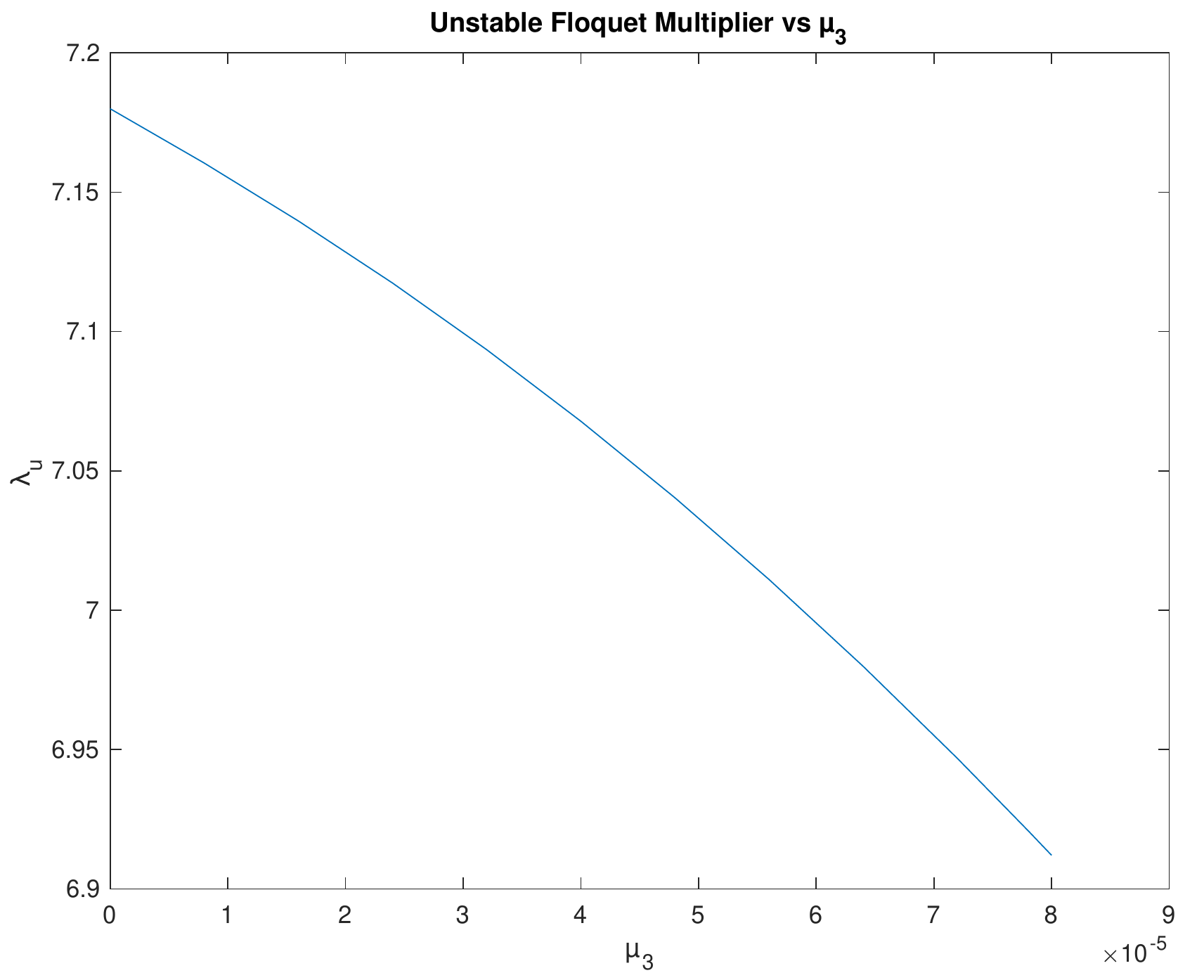}
\caption{ \label{fig:34lambdas} Unstable Floquet Multiplier vs $\mu_{3}$ for Jupiter-Europa 3:4 Continuation}
\end{centering}
\end{figure}

Zooming into the neighborhood of Europa in the $\mu_{3}$ continuation plot, it is visible that the invariant circle moves closer to Europa as $\mu_{3}$ increases, potentially to offset some of the gravitational pull of Ganymede in the other direction. In physical units, the closest approach to Europa decreases from 22052 km to 18721 km. Nevertheless, plotting the unstable Floquet multiplier $\lambda_{u}$ (recall Eq. \eqref{Lambdaform}) of the torus versus the parameter $\mu_{3}$ in Figure \ref{fig:34lambdas}, we can see that the instability of the torus actually decreases slightly as $\mu_{3}$ increases and the invariant circle moves towards Europa. 

\begin{figure}
\begin{centering}
\includegraphics[width=0.49\columnwidth]{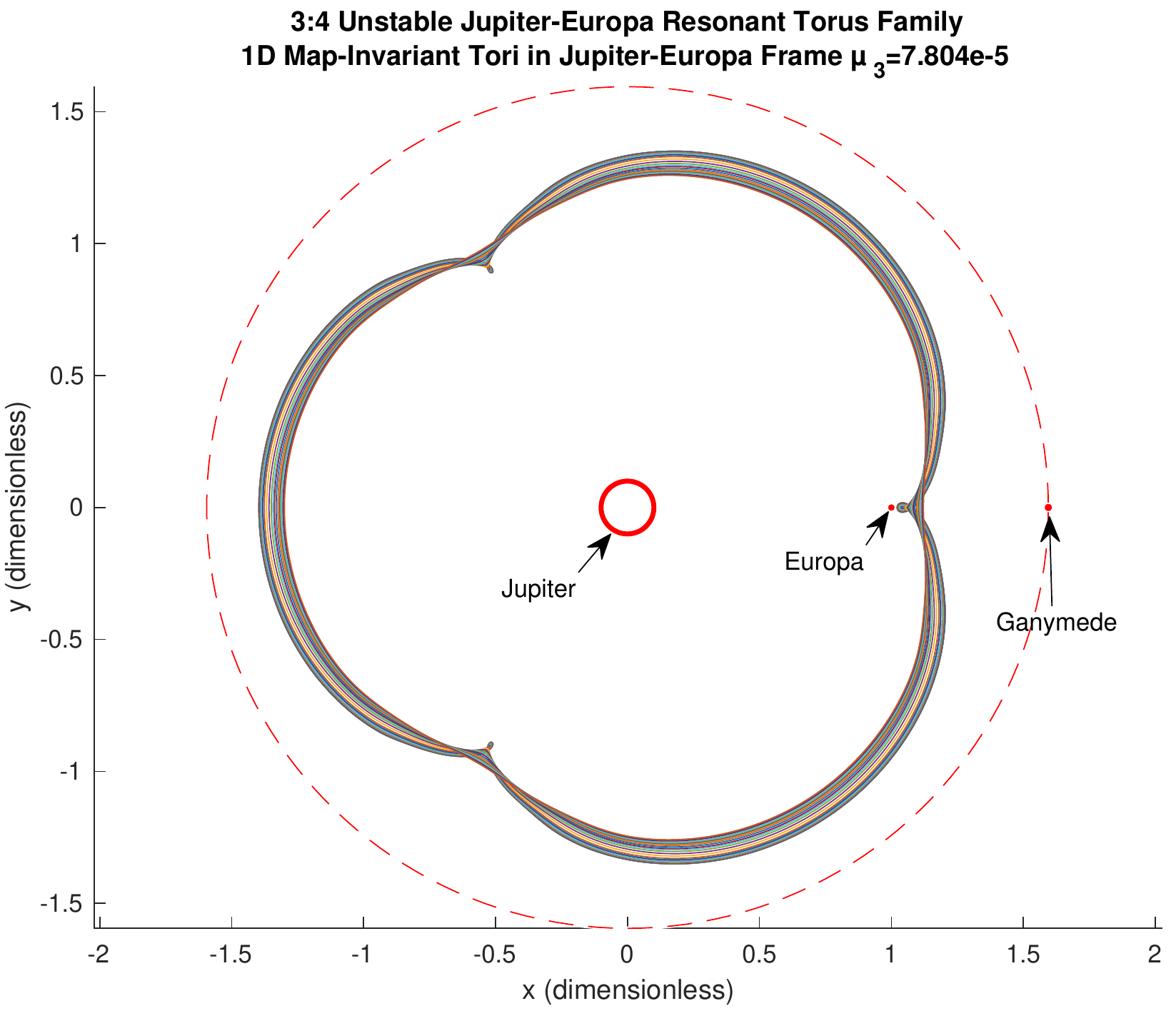}
\caption{ \label{fig:34omega} Continuation of Jupiter-Europa 3:4 Resonant CCR4BP Torus Family}
\end{centering}
\end{figure}

After computing a single Jupiter-Europa 3:4 unstable resonant CCR4BP torus for the physical value of Ganymede's mass ratio $\mu_{3}$, we continued this torus by $\omega$ to get a family of invariant circles for the stroboscopic map. Figure \ref{fig:34omega} shows the resulting tori computed in the Jupiter-Europa frame; the overall structure of the invariant circle family is similar to that of the Jupiter-Europa 3:4 unstable resonant PCRTBP periodic orbit family. 

\begin{figure}
\begin{centering}
\includegraphics[width=0.48\columnwidth]{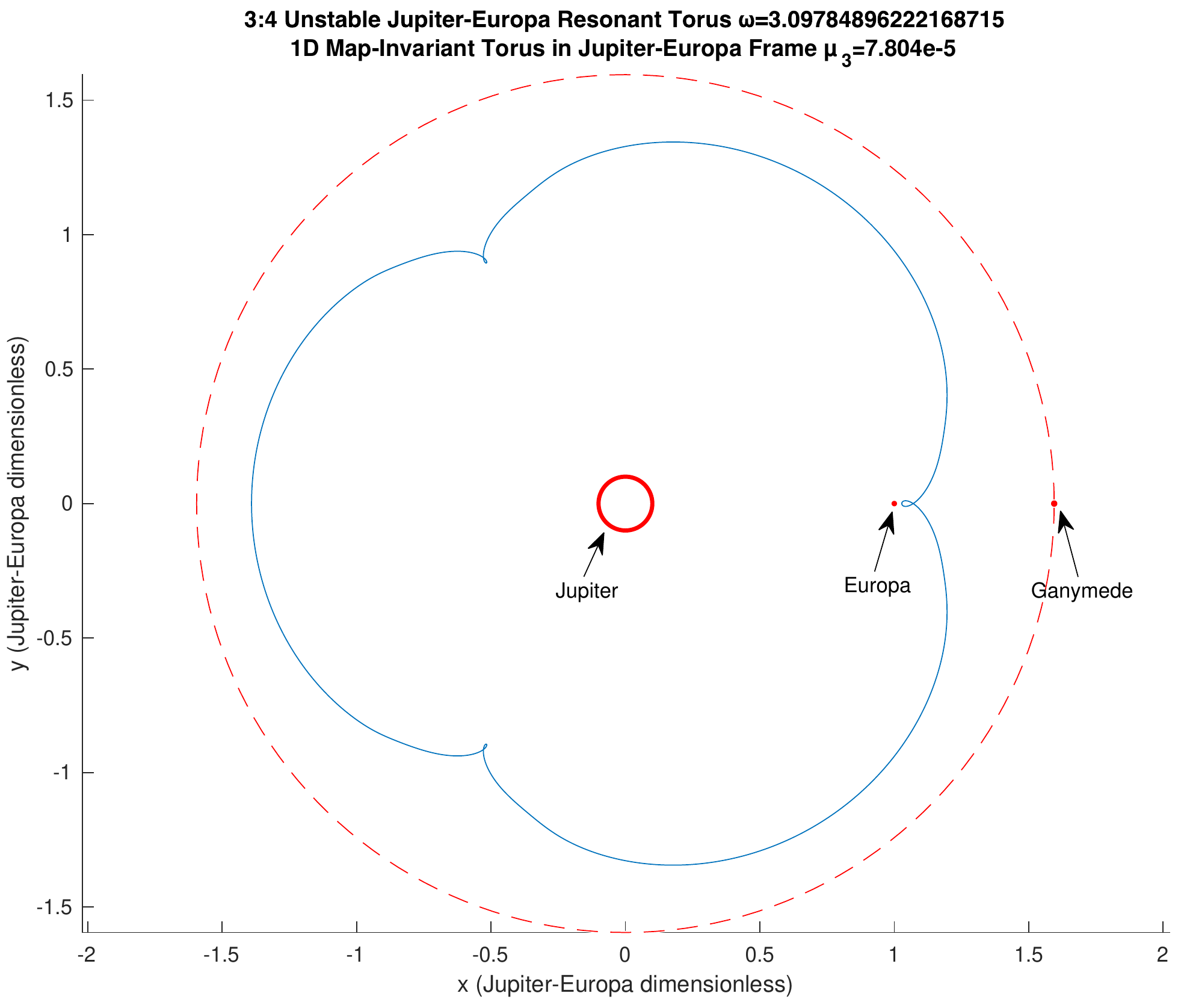}
\includegraphics[width=0.48\columnwidth]{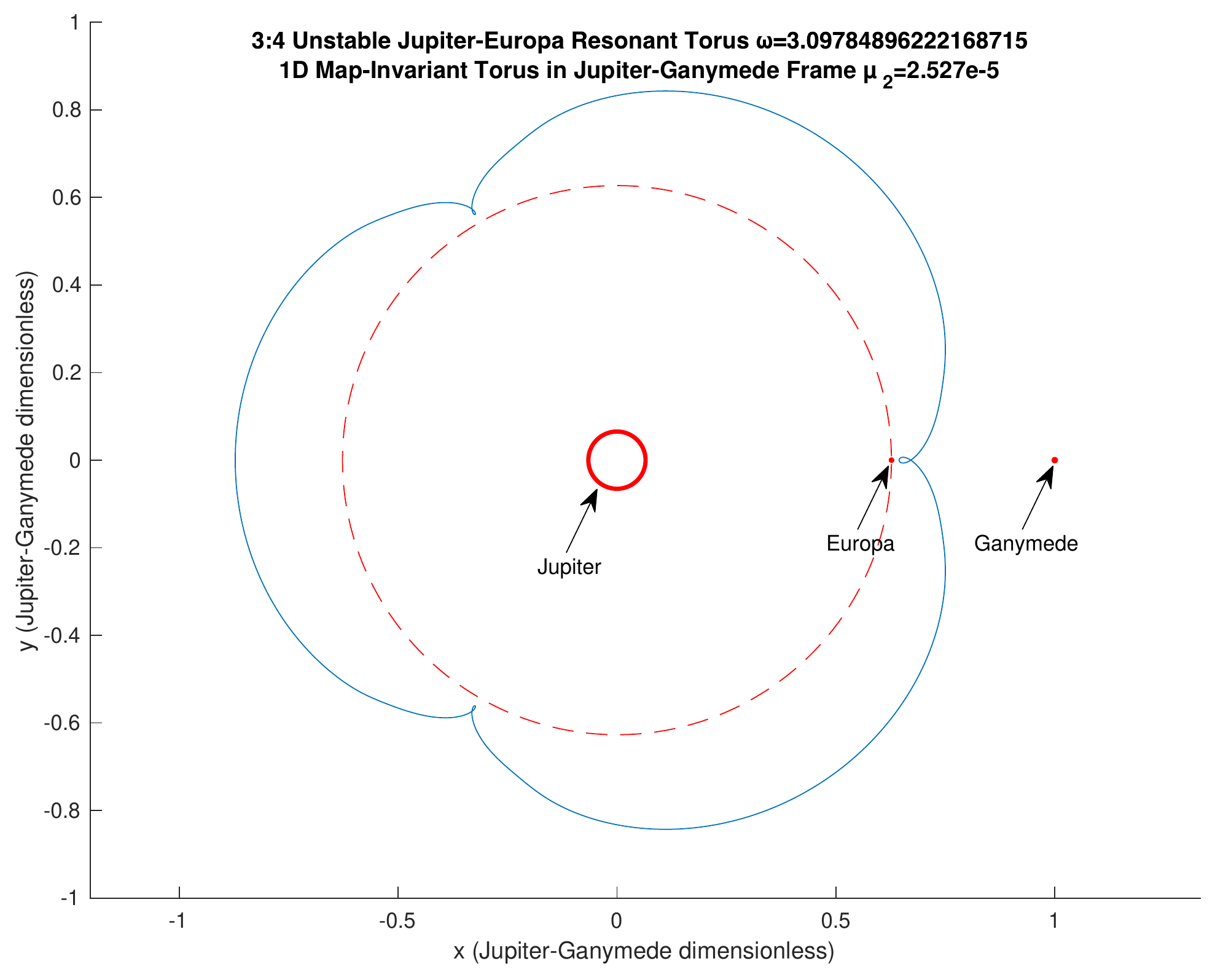}
\caption{ \label{fig:34maptori} 1D CCR4BP Stroboscopic Map-Invariant Torus for Jupiter-Europa 3:4 Resonance, (left) in the Jupiter-Europa frame, (right) in the Jupiter-Ganymede frame}
\end{centering}
\end{figure}

With Jupiter-Europa 3:4 unstable resonant CCR4BP invariant circles computed for the Jupiter-Europa frame stroboscopic map, we can transform the invariant circle states to the Jupiter-Ganymede rotating frame using the function $\Phi_{\theta_{3}}$ defined previously; note that one of the benefits of using $\theta_{3}=0$ as the stroboscopic map phase is that no rotation is required during this transformation between frames. As demonstrated earlier, the resulting set of states will be an invariant circle of the Jupiter-Ganymede frame CCR4BP stroboscopic map based at $\theta_{2}=0$. One such invariant circle is shown in both frames in Figure \ref{fig:34maptori}; on the left of the figure is the invariant circle in the Jupiter-Europa frame, while on the right is the same circle in the Jupiter-Ganymede frame. As is clearly visible, the shape of the invariant circle in position space is the same in either frame. 

\begin{figure}
\begin{centering}
\includegraphics[width=0.48\columnwidth]{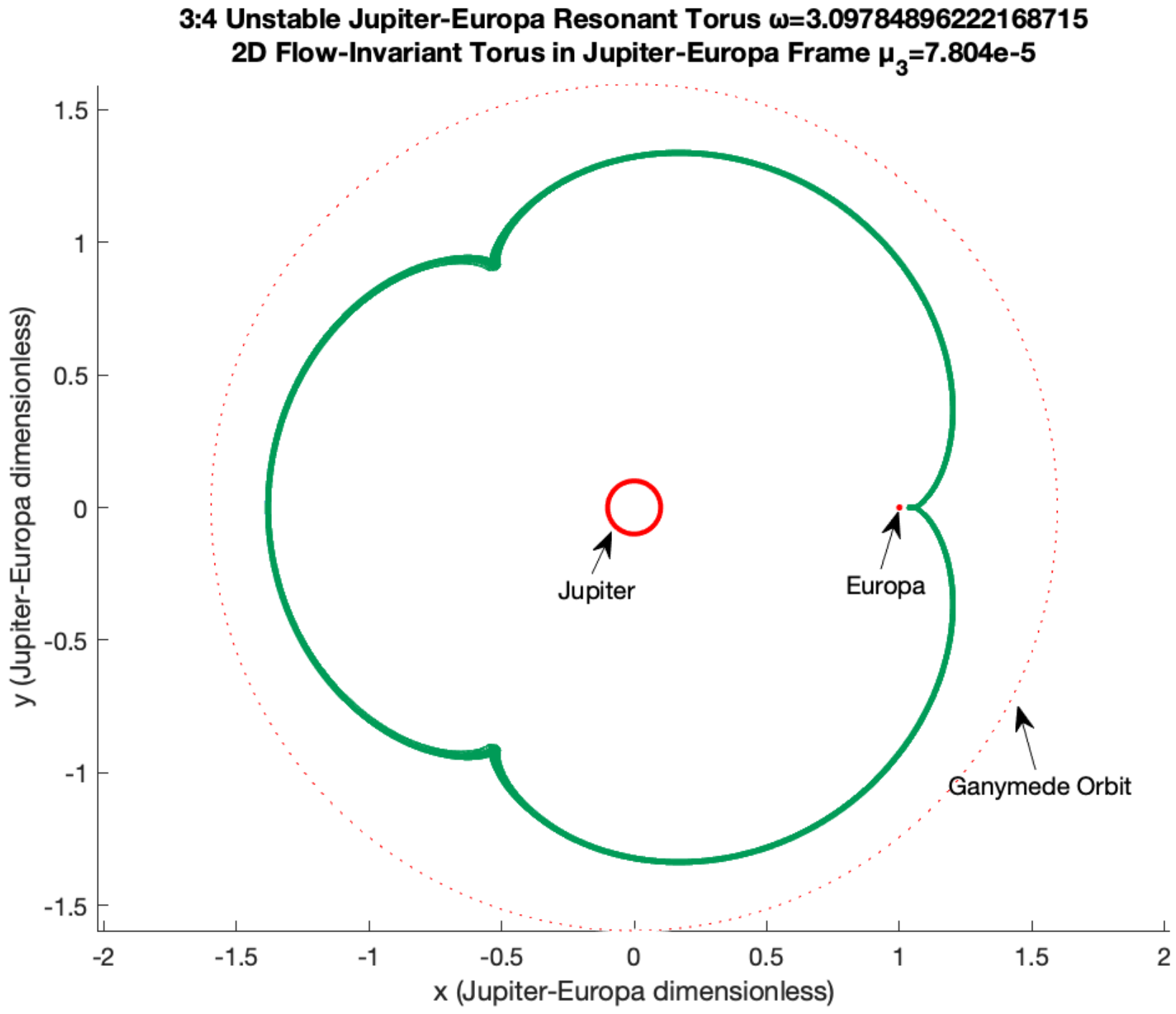}
\includegraphics[width=0.48\columnwidth]{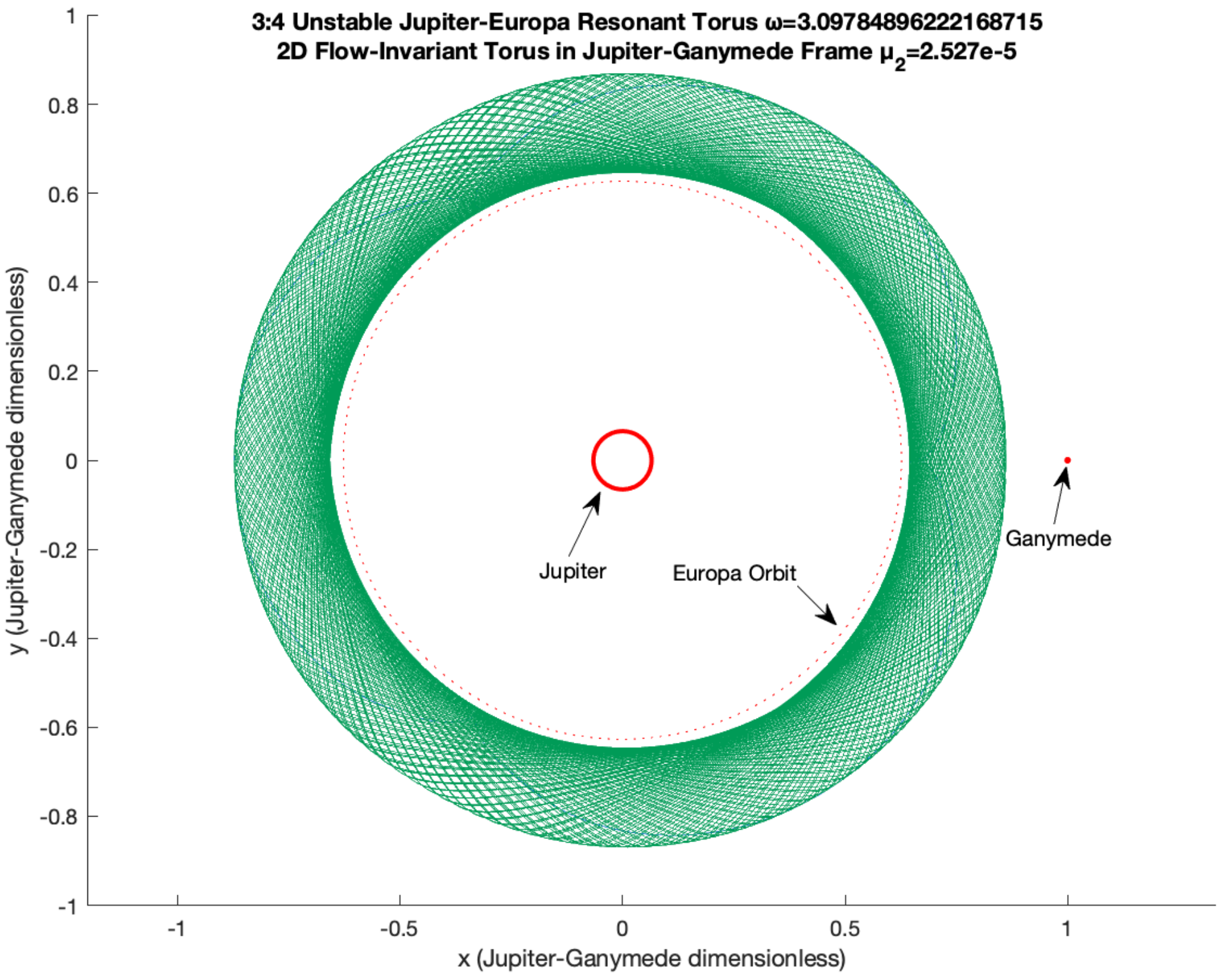}
\caption{\label{fig:34flowtori} 2D CCR4BP Flow-Invariant Torus for Jupiter-Europa 3:4 Resonance, (left) in the Jupiter-Europa frame, (right) in the Jupiter-Ganymede frame}
\end{centering}
\end{figure}

We can also convert this 1D CCR4BP stroboscopic map invariant circle to a 2D invariant torus of the continuous time equations of motion in either frame; this is done by integrating the circle states by their respective CCR4BP equations of motion. The resulting 2D torus corresponding to the invariant circle of Figure \ref{fig:34maptori} is plotted in both frames in Figure \ref{fig:34flowtori};  it is clear that the 2D flow-invariant torus looks very different depending on the frame used, even though the map-invariant circle was similarly shaped in both frames. Indeed, in the Jupiter-Europa frame, the continuous-time trajectories closely follow the shape of the invariant circle, whereas in the Jupiter-Ganymede frame, a much larger portion of the position space is explored.

\begin{figure}
\begin{centering}
\includegraphics[width=0.49\columnwidth]{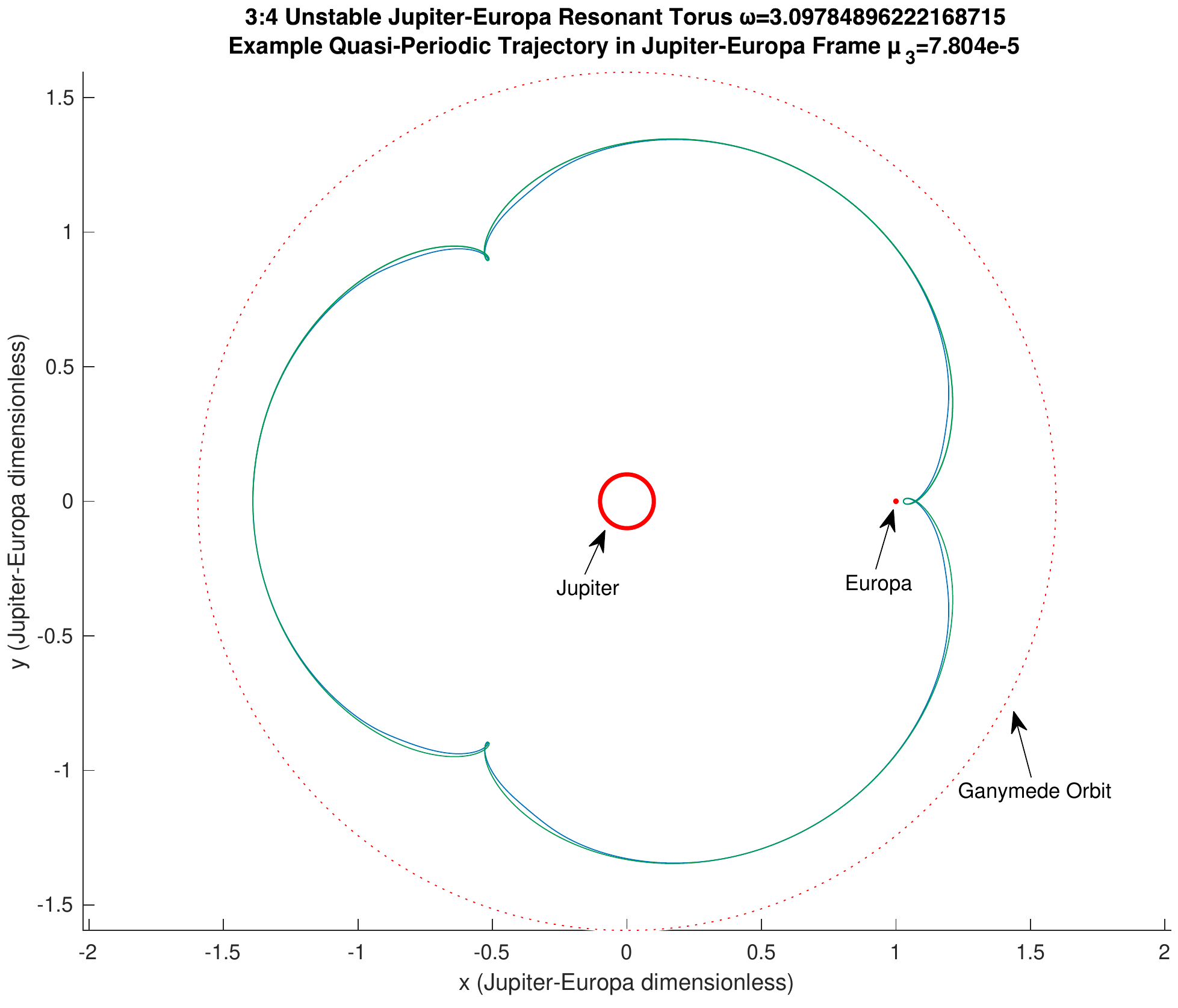}
\includegraphics[width=0.49\columnwidth]{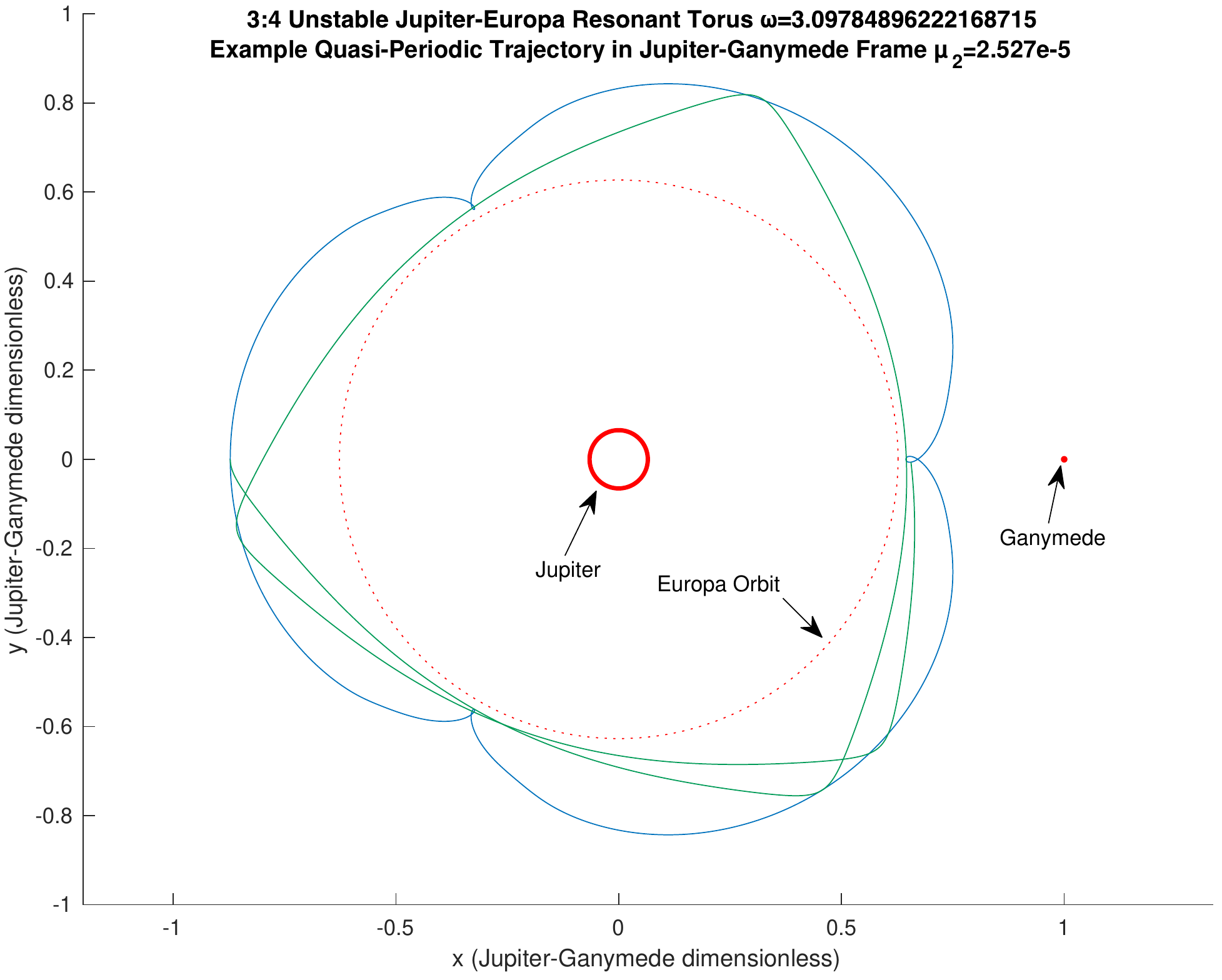}
\caption{ \label{fig:34toritraj} Example Jupiter-Europa 3:4 CCR4BP Quasi-Periodic Trajectory, (left) in the Jupiter-Europa frame, (right) in the Jupiter-Ganymede frame}
\end{centering}
\end{figure}

To better visualize what is happening, we numerically integrated a single state from the invariant circle to get a continuous-time trajectory lying on the torus of Figure \ref{fig:34flowtori}; this was then plotted in each frame, as is shown in Figure \ref{fig:34toritraj}. The invariant circle is plotted in blue, while the integrated trajectory is in green; the trajectory initial condition was chosen to be the intersection of the invariant circle with the negative $x$-axis. The integration time in the Jupiter-Europa frame was $\frac{6\pi}{|\Omega_3-1|}$, while in the Jupiter-Ganymede frame it was $\frac{6\pi}{|\Omega_2-1|}$; these are the same physical amount of time and are related to each other through the change of time unit between frames. Comparing the trajectory plotted in the two frames, it is visible that the trajectory's apoapsis passes seem to rotate clockwise in the Jupiter-Ganymede frame over time; by contrast, in the Jupiter-Europa frame, the apoapsis passes remain in the vicinity of the same three locations over time. 

\subsection{3:2 Jupiter-Ganymede Resonant CCR4BP Tori}

In Figure \ref{fig:34toritraj}, plotting the 3:4 Jupiter-Europa resonant CCR4BP torus trajectory in the Jupiter-Ganymede frame gave a trajectory somewhat reminiscent of a 3:2 Jupiter-Ganymede PCRTBP resonant periodic orbit. Hence, we next computed 3:2 Jupiter-Ganymede unstable resonant CCR4BP tori. During our first attempts at continuing 3:2 Jupiter-Ganymede unstable PCRTBP periodic orbits into invariant tori of the CCR4BP (we used the Jupiter-Ganymede frame for this), we found our quasi-Newton method failing to converge for even extremely small $\bar \mu_{2}$ values such as $10^{-12}$. 

\begin{figure}
\begin{centering}
\includegraphics[width=0.49\columnwidth]{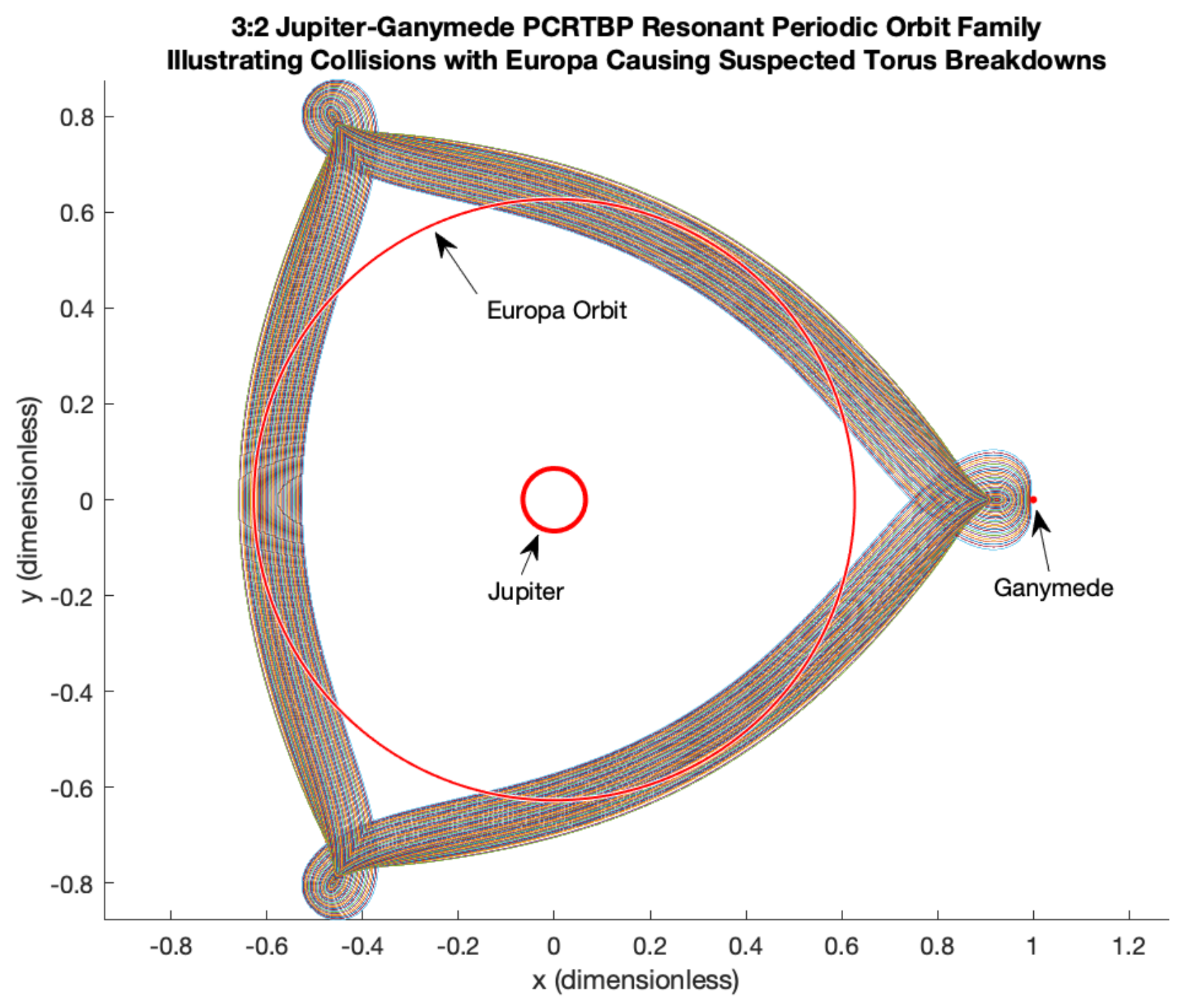}
\caption{ \label{fig:32collisions} 3:2 Jupiter-Ganymede Periodic Orbit Collisions with Europa}
\end{centering}
\end{figure}

Upon further investigation, the reason for this was found to be that the addition of Europa to the model introduces singularities in the equations of motion at points belonging to many of the 3:2 Jupiter-Ganymede periodic orbits. In Figure \ref{fig:32collisions}, we plot part of the family of unstable 3:2 Jupiter-Ganymede PCRTBP resonant orbits, as well as a red circle representing the orbit of Europa in the CCR4BP. It is clearly visible that the Europa orbit intersects all of the 3:2 Jupiter-Ganymede orbits which pass closest to Ganymede. Hence, no matter how small we take $\mu_{2}$, unless $\bar \mu_{2}=0$, we are introducing an infinite perturbation of the equations of motion at states which are included in those periodic orbits; we suspect this destroys any potential invariant tori which would have existed near those periodic orbits for a finite perturbation. We have not investigated whether such tori might exist if one uses regularized equations of motion to remove the singularity at Europa, however. 

\begin{figure}
\begin{centering}
\includegraphics[width=0.49\columnwidth]{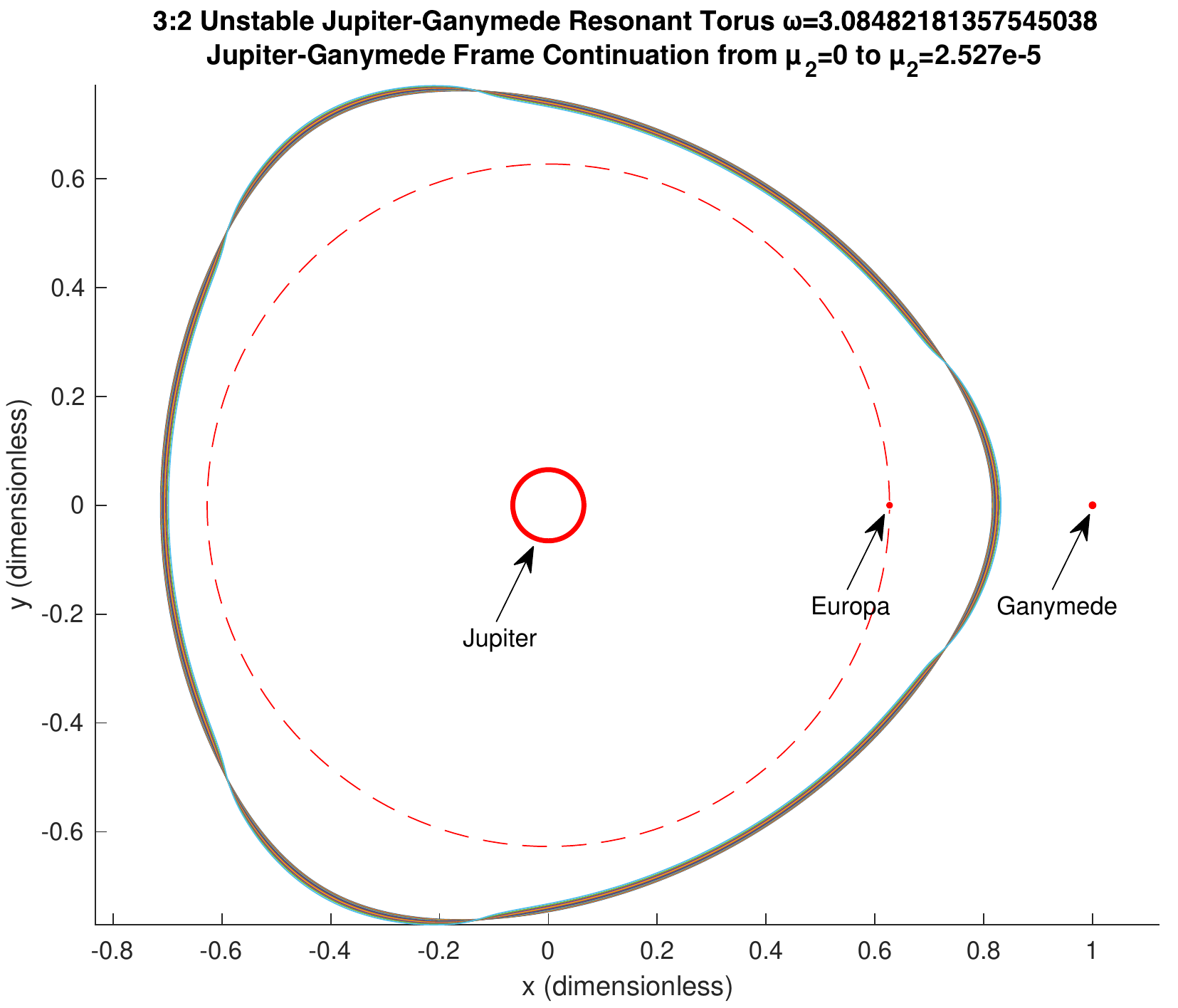}
\includegraphics[width=0.49\columnwidth]{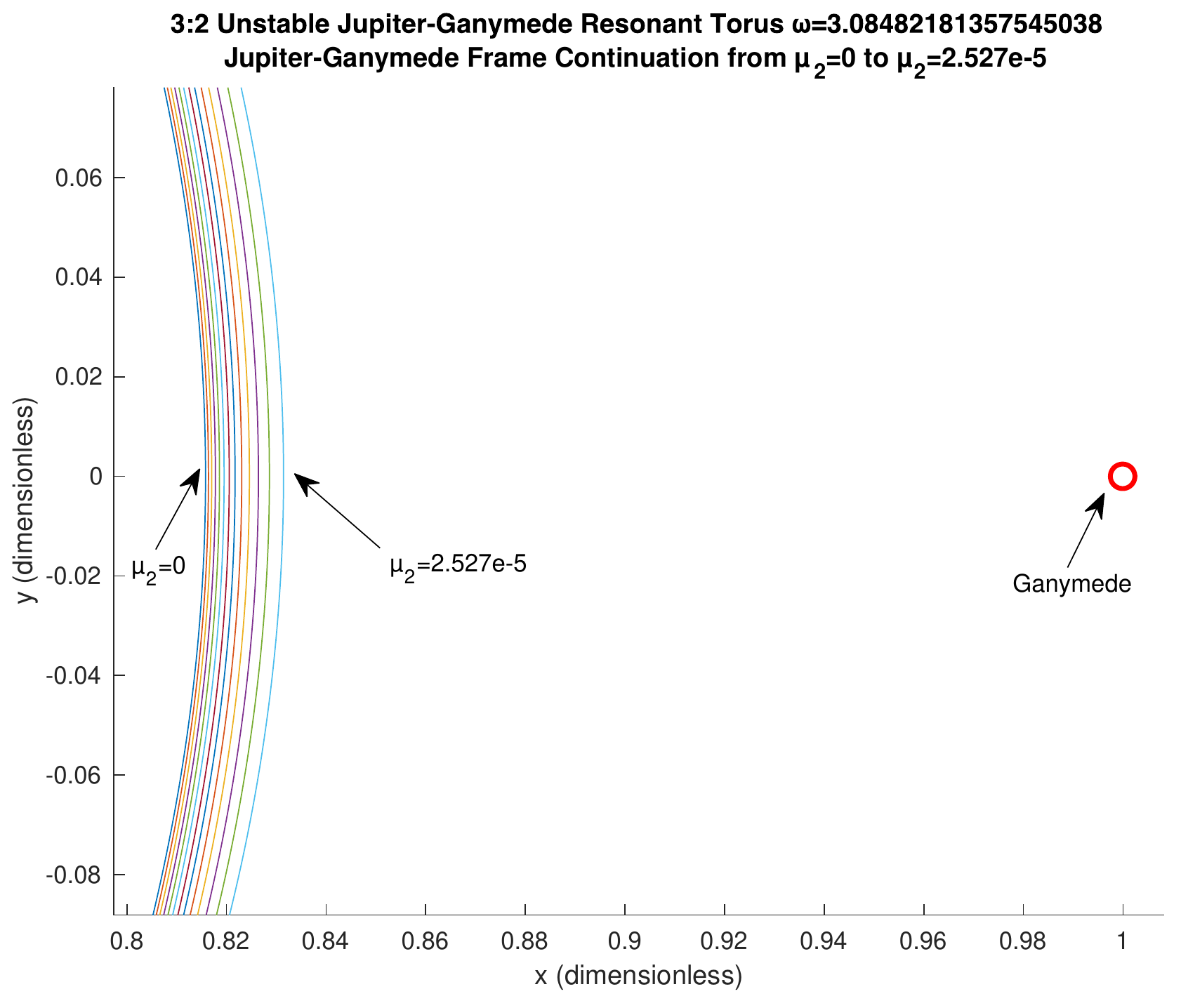}
\caption{ \label{fig:32eps} Continuation of Jupiter-Ganymede 3:2 Resonant CCR4BP Torus By $\mu_{2}$}
\end{centering}
\end{figure}

While all of the most unstable 3:2 Jupiter-Ganymede PCRTBP unstable resonant periodic orbits intersect Europa's orbit, some of the less unstable ones do not; some of these are visible as the outermost orbits in Figure \ref{fig:32collisions}, having higher perijoves and lower apojoves. We were able to continue one such orbit by $\bar \mu_{2}$ to the physical value of $\bar \mu_{2}=2.5265115494603433 \times 10^{-5}$ for Europa in the Jupiter-Ganymede frame planar CCR4BP. The continuation step size used was $\Delta \bar \mu_{2} = 2.5 \times 10^{-7}$. Some of the computed intermediate stroboscopic map invariant circles are shown in Figure \ref{fig:32eps},  with a zoomed-in plot near Ganymede given on the right. 

\begin{figure}
\begin{centering}
\includegraphics[width=0.49\columnwidth]{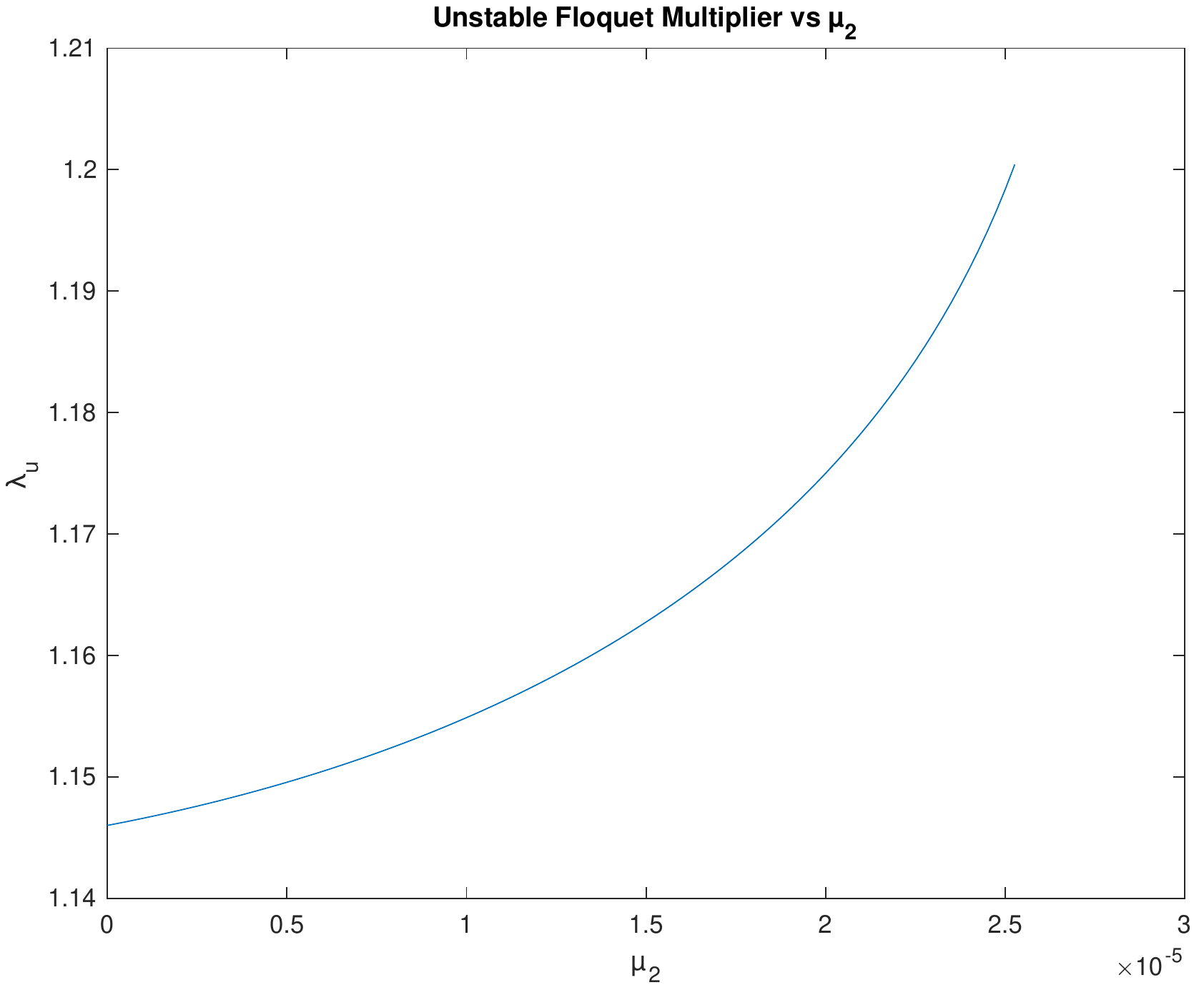}
\caption{ \label{fig:32lambdas} Unstable Floquet Multiplier vs $\bar \mu_{2}$ for Jupiter-Ganymede 3:2 Continuation}
\end{centering}
\end{figure}

Just as was seen in the Jupiter-Europa 3:4 resonant CCR4BP torus computation, we see the invariant circles moving closer to Ganymede as $\bar \mu_{2}$ is increased. However, this time, plotting the torus' unstable Floquet multiplier $\lambda_{u}$ versus the parameter $\bar \mu_{2}$ in this case shows the opposite trend compared to the Jupiter-Europa 3:4 resonance. In Figure \ref{fig:32lambdas}, it can be seen that as $\bar \mu_{2}$ increases and the circle's closest approach moves closer and closer to Ganymede, the torus grows more unstable. However, it should be noted that the value of $\lambda_{u}$ is much smaller for this 3:2 Jupiter-Ganymede resonant torus than it was for the 3:4 Jupiter-Europa torus.

\subsection{7:5 Jupiter-Ganymede Resonant CCR4BP Tori}

\begin{figure}
\begin{centering}
\includegraphics[width=0.49\columnwidth]{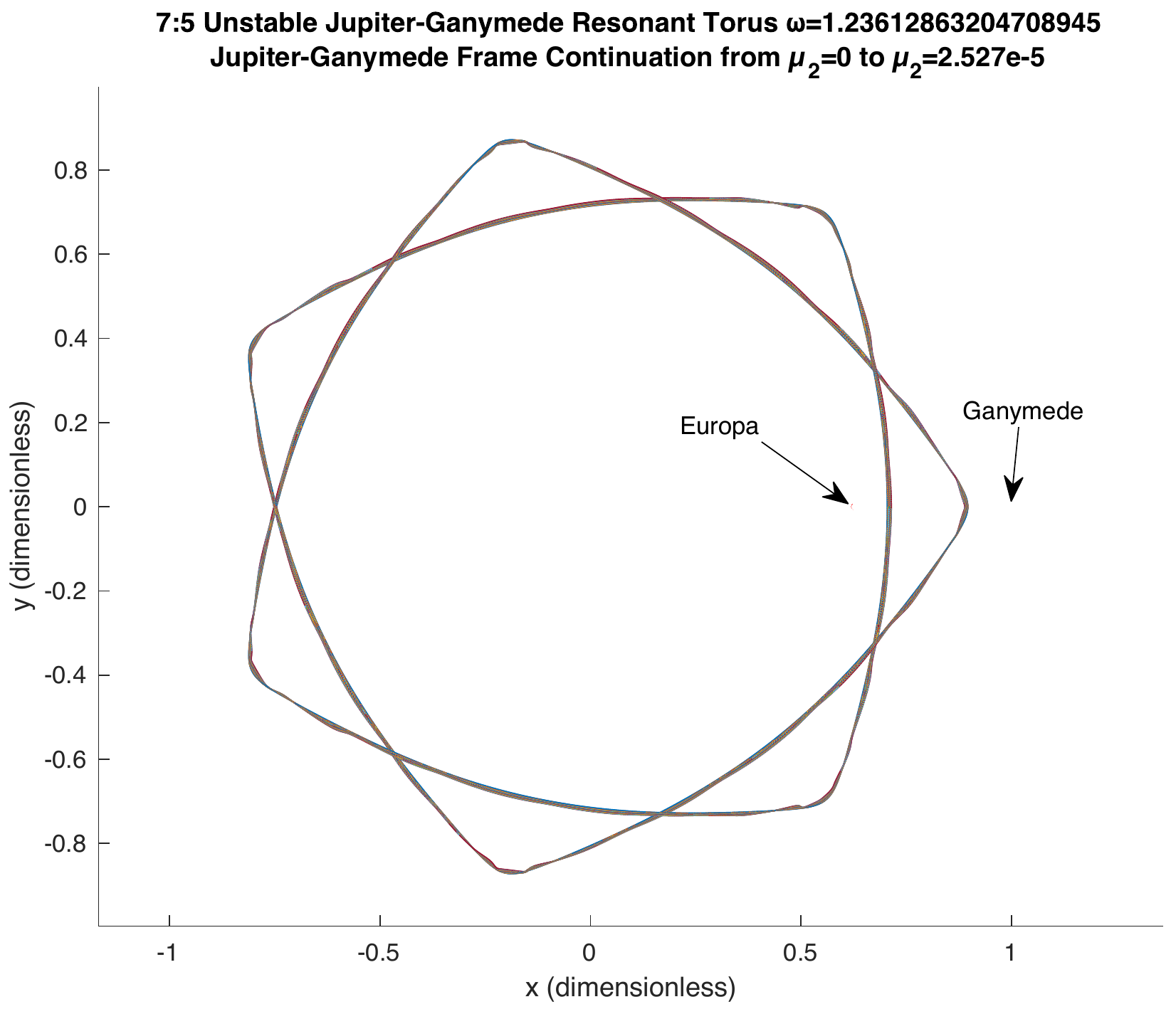}
\includegraphics[width=0.49\columnwidth]{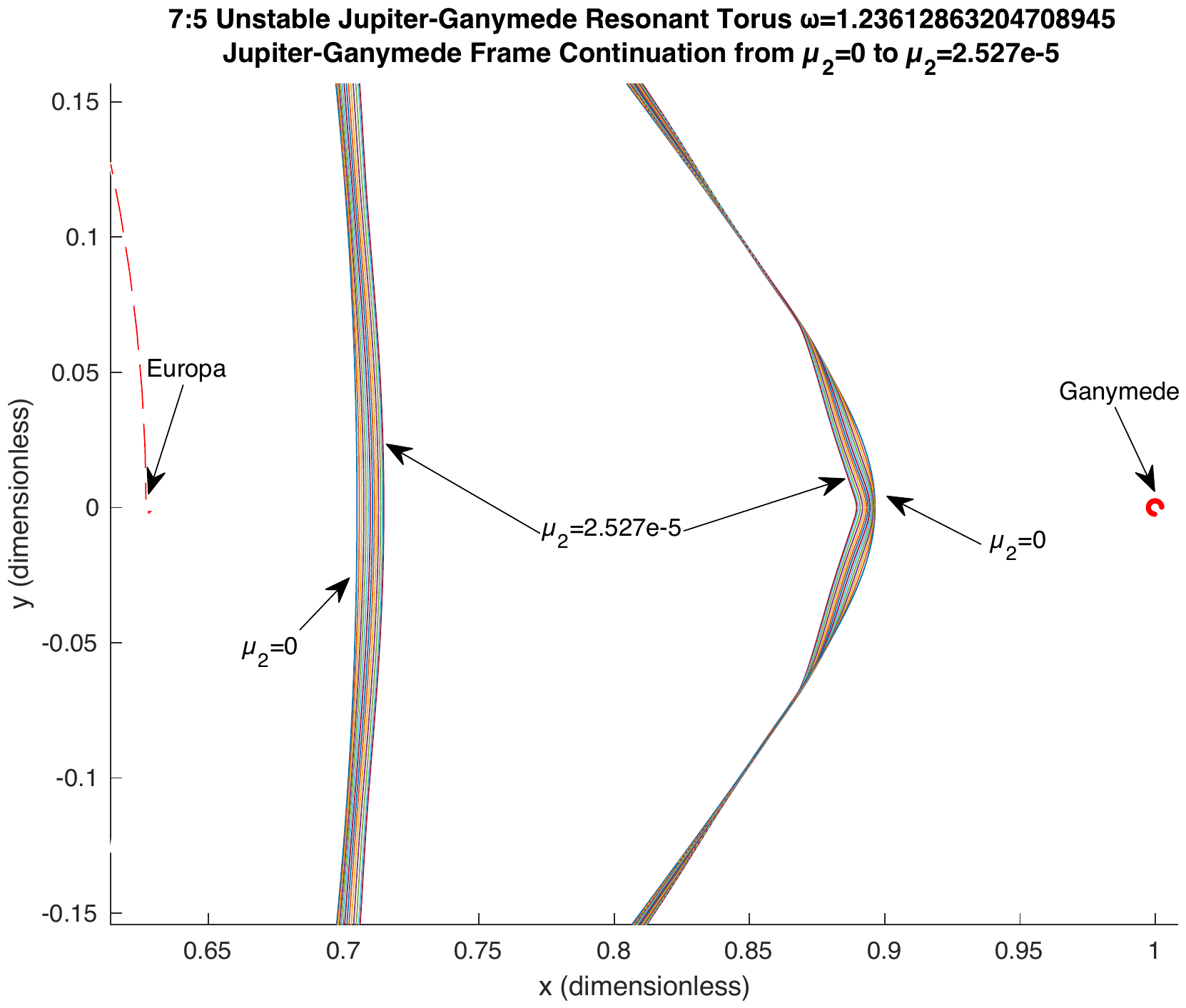}
\caption{ \label{fig:75eps} Continuation of Jupiter-Ganymede 7:2 Resonant CCR4BP Torus By $\mu_{2}$}
\end{centering}
\end{figure}

Another resonance of interest for Jupiter-Europa-Ganymede tour design is the Jupiter-Ganymede 7:5 resonance\cite{Anderson2021a}. Hence, we computed the corresponding tori in the CCR4BP model. Figure \ref{fig:75eps} shows the continuation of one such invariant circle by $\bar \mu_{2}$ to the physical value for Europa in the Jupiter-Ganymede frame. As $\bar \mu_{2}$ went from zero to its physical value, the number of Fourier modes required to sufficiently accurately represent the torus increased significantly, from 512 to 2048. 

\begin{figure}
\begin{centering}
\includegraphics[width=0.48\columnwidth]{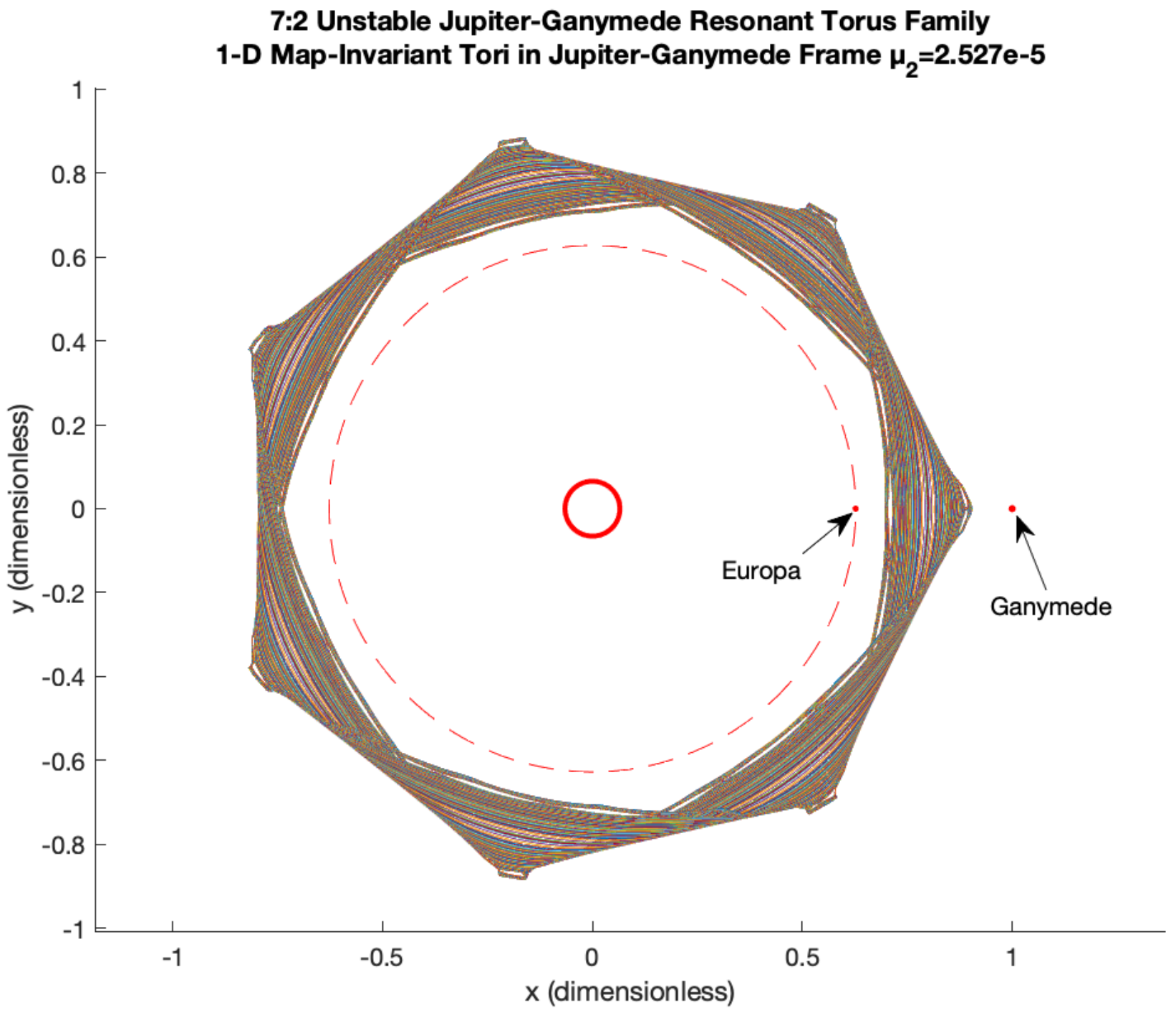}
\caption{ \label{fig:75omega} Family of Jupiter-Ganymede 7:2 Resonant CCR4BP Tori}
\end{centering}
\end{figure}

\begin{figure}
\begin{centering}
\includegraphics[width=0.48\columnwidth]{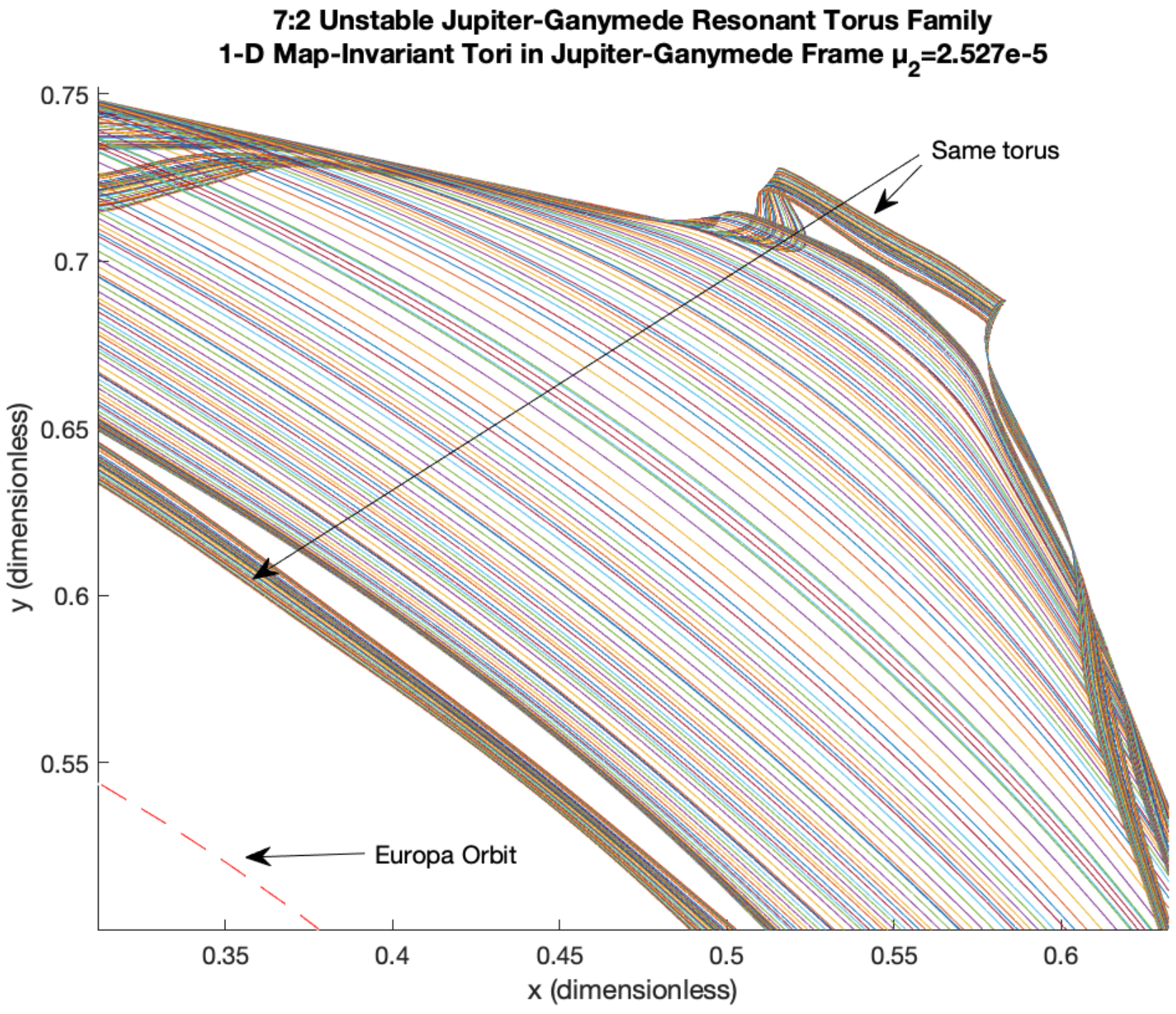}
\includegraphics[width=0.48\columnwidth]{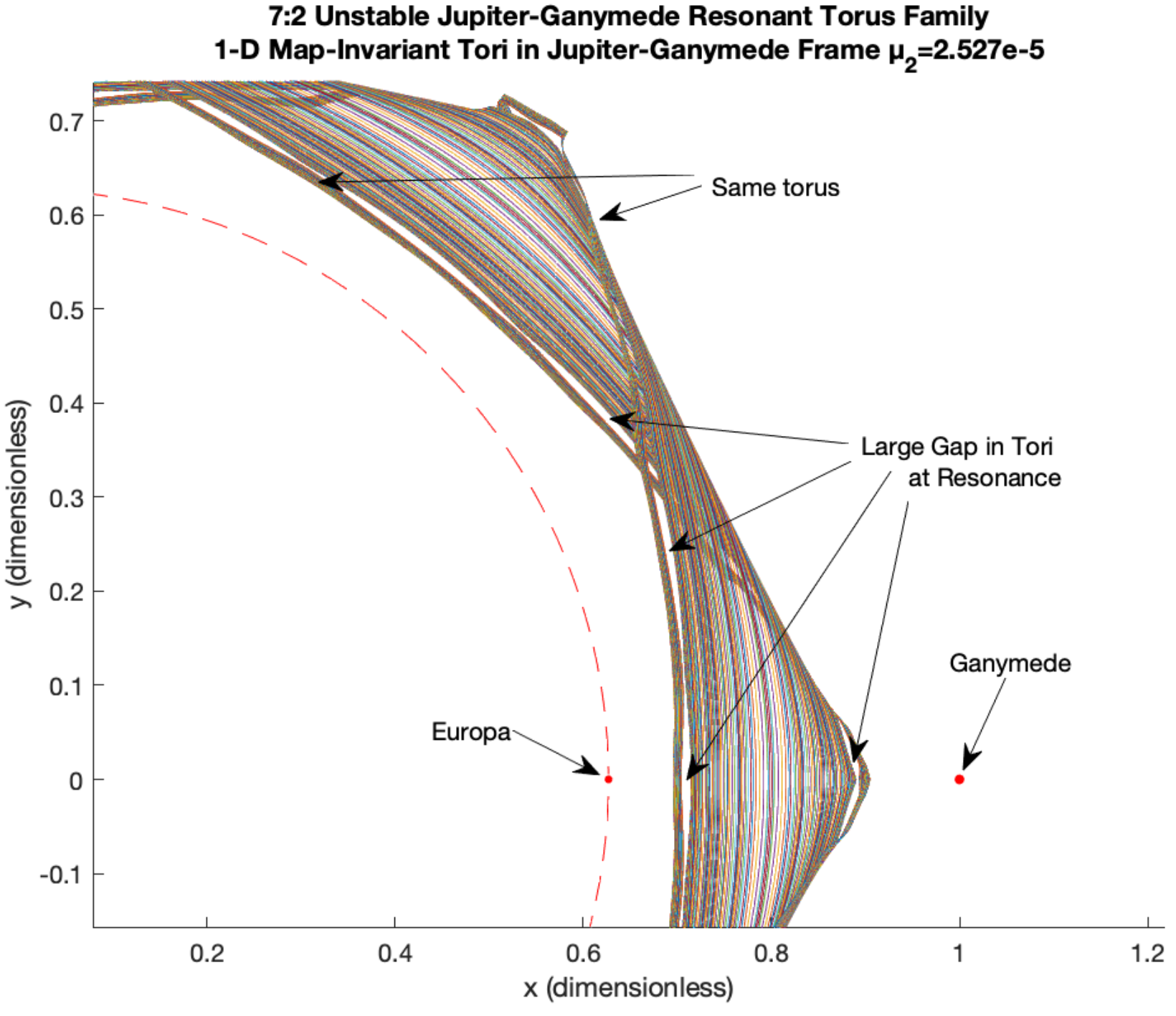}
\caption{ \label{fig:75omegazoomed} Continuation of Jupiter-Ganymede 3:2 Resonant CCR4BP Torus By $\mu_{2}$}
\end{centering}
\end{figure}

Using a combination of continuation by $\bar \mu_{2}$ and $\omega$, we next computed a family of 7:5 Jupiter-Ganymede resonant tori in the physical Jupiter-Europa-Ganymede CCR4BP; this is plotted in Figure \ref{fig:75omega}. Zooming into the family as shown in Figure \ref{fig:75omegazoomed}, we see some interesting phenomena. First of all, it is clear that certain invariant circles see sharp bends at their apojoves; as indicated in the left plot, these apojoves belong to the same invariant circles which make closer passes of Europa's orbit (note that the 7:5 invariant circles loop twice around Jupiter in position space before closing). 

Second of all, as indicated on the right plot of Figure \ref{fig:75omegazoomed}, we see a significant gap between the invariant circles we were able to compute. The reason for this is that in this region, the rotation number of topologically similar tori (had they existed) would have been nearly rational; this was verified by computing the continued fraction expansion of the average rotation number of the two tori on either side of the gap. Equivalently, the period of the perturbation from Europa is nearly resonant with the period of the corresponding Jupiter-Ganymede PCRTBP periodic orbits. This causes torus breakdown. Indeed, the ``bulge" shape of the gap is immediately reminiscent of the gaps between invariant circles which appear at rational rotation numbers in the commonly-studied Chirikov standard map\cite{sander2020}, as the perturbation strength is increased. 

\section{Relations Between Jupiter-Europa and Jupiter-Ganymede Resonant Tori} \label{relationsSection}

\begin{figure}
\begin{centering}
\includegraphics[width=0.49\columnwidth]{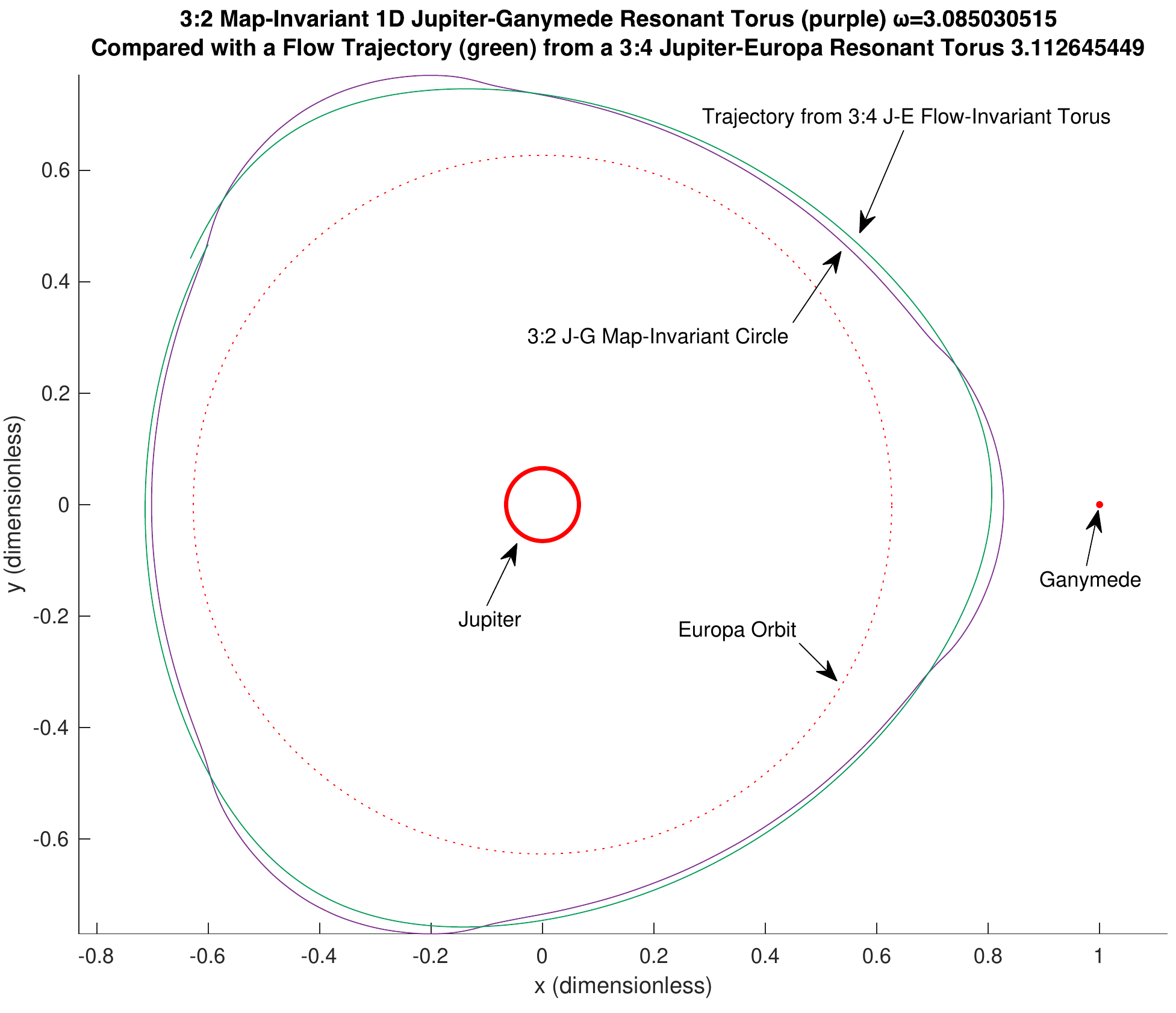}
\caption{ \label{fig:3432traj} Comparison of a Flow Trajectory from a Jupiter-Europa 3:4 Torus to a Jupiter-Ganymede 3:2 Invariant Circle in the Jupiter-Ganymede Rotating Frame}
\end{centering}
\end{figure}

As mentioned earlier, the 3:4 Jupiter-Europa resonant CCR4BP torus trajectory shown in Figure \ref{fig:34toritraj} in the Jupiter-Ganymede frame gave a trajectory of a similar shape to a 3:2 Jupiter-Ganymede PCRTBP resonant periodic orbit. Indeed, we actually found a continuous-time trajectory which starts on one of the 3:4 Jupiter-Europa family of CCR4BP tori (see Figure \ref{fig:34omega}), which also for some time closely follows another 3:2 Jupiter-Ganymede CCR4BP invariant circle (as well as the continuous-time trajectories of the corresponding 2D flow-invariant torus). The aforementioned 3:4 Jupiter-Europa trajectory and the 3:2 Jupiter-Ganymede map-invariant circle are plotted in green and blue, respectively, in a Jupiter-Ganymede rotating frame in Figure \ref{fig:3432traj}. 

The reason for this observed similarity in the shape of Jupiter-Europa 3:4 and Jupiter-Ganymede 3:2 resonant CCR4BP torus trajectories is the Laplace resonance\cite{Barnes2011}. The Laplace resonance refers to the 4:2:1 resonance between the orbits of Io, Europa, and Ganymede around Jupiter; in other words, Io makes approximately four revolutions around Jupiter in the same time that Europa makes two and Ganymede makes one. Hence, if the spacecraft makes approximately $m$ revolutions around Jupiter in the same time that Ganymede makes $n$ revolutions (an $m$:$n$ resonance), this is also approximately the same amount of time in which Europa makes $2n$ revolutions. This implies that the spacecraft is close to an $m$:$2n$ resonance with Europa if it is in an $m$:$n$ resonance with Ganymede. Of course, the previously observed similarity between Jupiter-Ganymede 3:2 and Jupiter-Europa 3:4 trajectories is just the case of $m=3$ and $n=2$. 

Given the potential for a relationship between Jupiter-Ganymede $m$:$n$ and Jupiter-Europa $m$:$2n$ resonances in the planar CCR4BP, we carried out an investigation to see what happens if we take a Jupiter-Europa 3:4 resonant torus in the CCR4BP with full Jupiter, Europa, and Ganymede masses, and continue this to the Jupiter-Ganymede PCRTBP by decreasing the Europa mass to zero. Hence, the overall continuation path taken is of continuing a Jupiter-Europa 3:4 resonant periodic orbit by $\mu_{3}$ from the Jupiter-Europa PCRTBP to the CCR4BP, and then (after changing frames) continuing the result by $\bar \mu_{2}$ from the CCR4BP to the Jupiter-Ganymede PCRTBP. 

We started with the invariant circle from Figure \ref{fig:34omega} with $\omega=3.111756$, and after transforming to the Jupiter-Ganymede rotating frame, started the decreasing-$\bar \mu_{2}$ continuation  from the physical value $\bar \mu_{2} = 2.5265 \times 10^{-5}$. One thing to note is that in the PCRTBP, the range of possible periods for the family of 3:4 resonant periodic orbits varies as the mass ratio changes (indeed, for zero mass ratio in the PCRTBP, the period of all 3:4 resonant orbits must be exactly $8\pi$); this change of periods depending on Europa mass corresponds to changes in the range of possible rotation numbers of our tori as $\mu_{2}$ decreases. In the limit of zero Europa mass, we expect the torus frequency $\Omega_1$ from Eq. \eqref{torusXTraj} corresponding to the original Jupiter-Europa 3:4 periodic orbit to tend to $\Omega_1=\frac{2\pi}{8\pi}\frac{\Omega_2}{\Omega_3} = 0.503493$ in Jupiter-Ganymede frame time units. Thus, $\omega$ is expected to tend towards $\omega=\frac{2\pi\Omega_1}{\Omega_2-1}=3.119948$ for all 3:4 tori continued to $\bar \mu_{2}=0$. Since we started the continuation at $\omega = 3.111756$, we had to interrupt our $\bar \mu_{2}$ continuation in the middle to continue $\omega$ to larger values as well. 

\begin{figure}
\begin{centering}
\includegraphics[width=0.48\columnwidth]{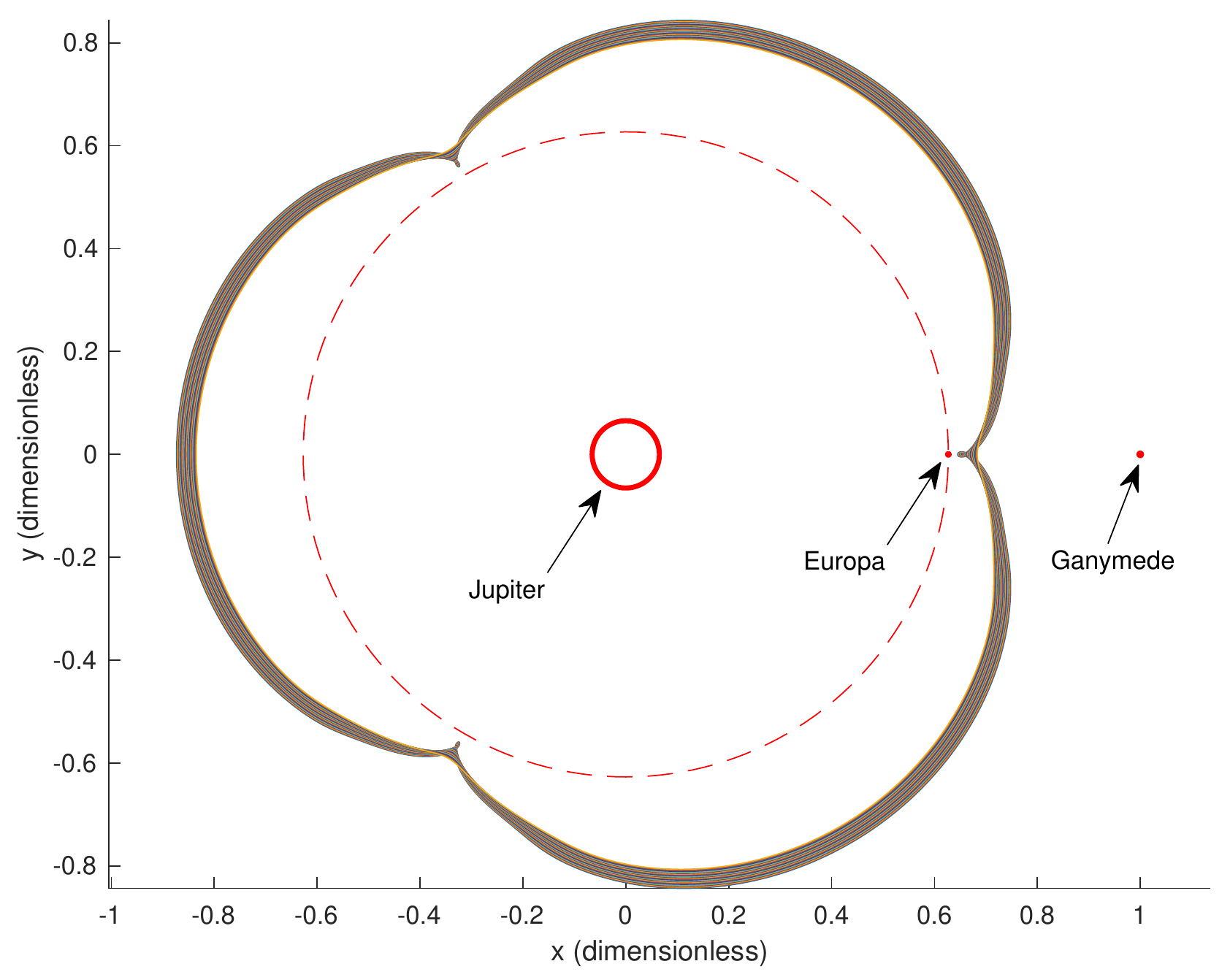}
\includegraphics[width=0.48\columnwidth]{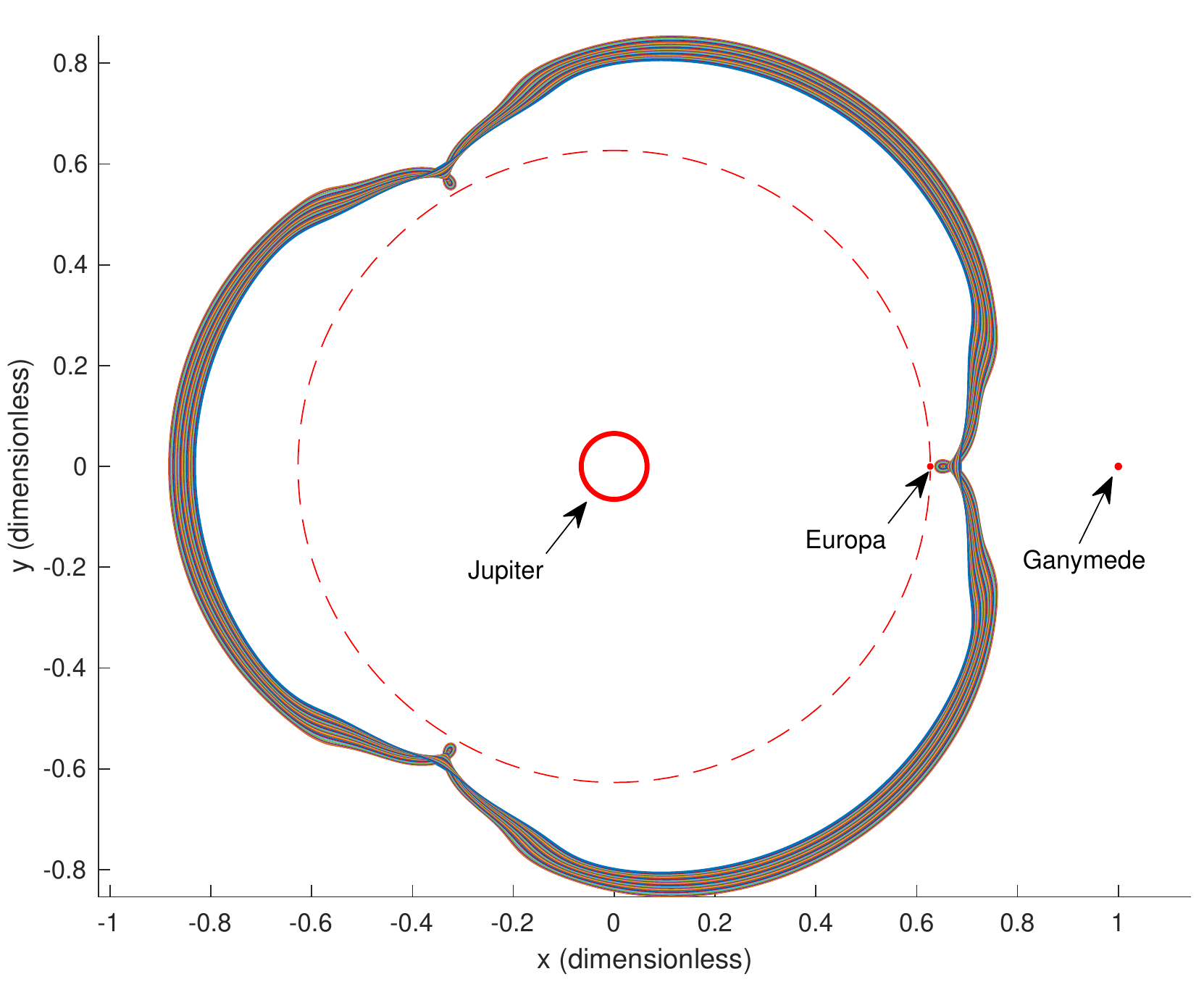}
\caption{ \label{fig:34epsdown} Continuation by $\bar \mu_{2}$ of Jupiter-Europa 3:4 resonant invariant circle in Jupiter-Ganymede frame from (left) $\bar \mu_{2} = 2.5265 \times 10^{-5}$ to $1.0015 \times 10^{-5}$ with $\omega = 3.111756$; (right) $\bar \mu_{2} = 1.0015 \times 10^{-5}$ to $2.506370 \times 10^{-6}$ with $\omega = 3.116809$
}
\end{centering}
\end{figure} 
\begin{figure}
\begin{centering}
\includegraphics[width=0.48\columnwidth]{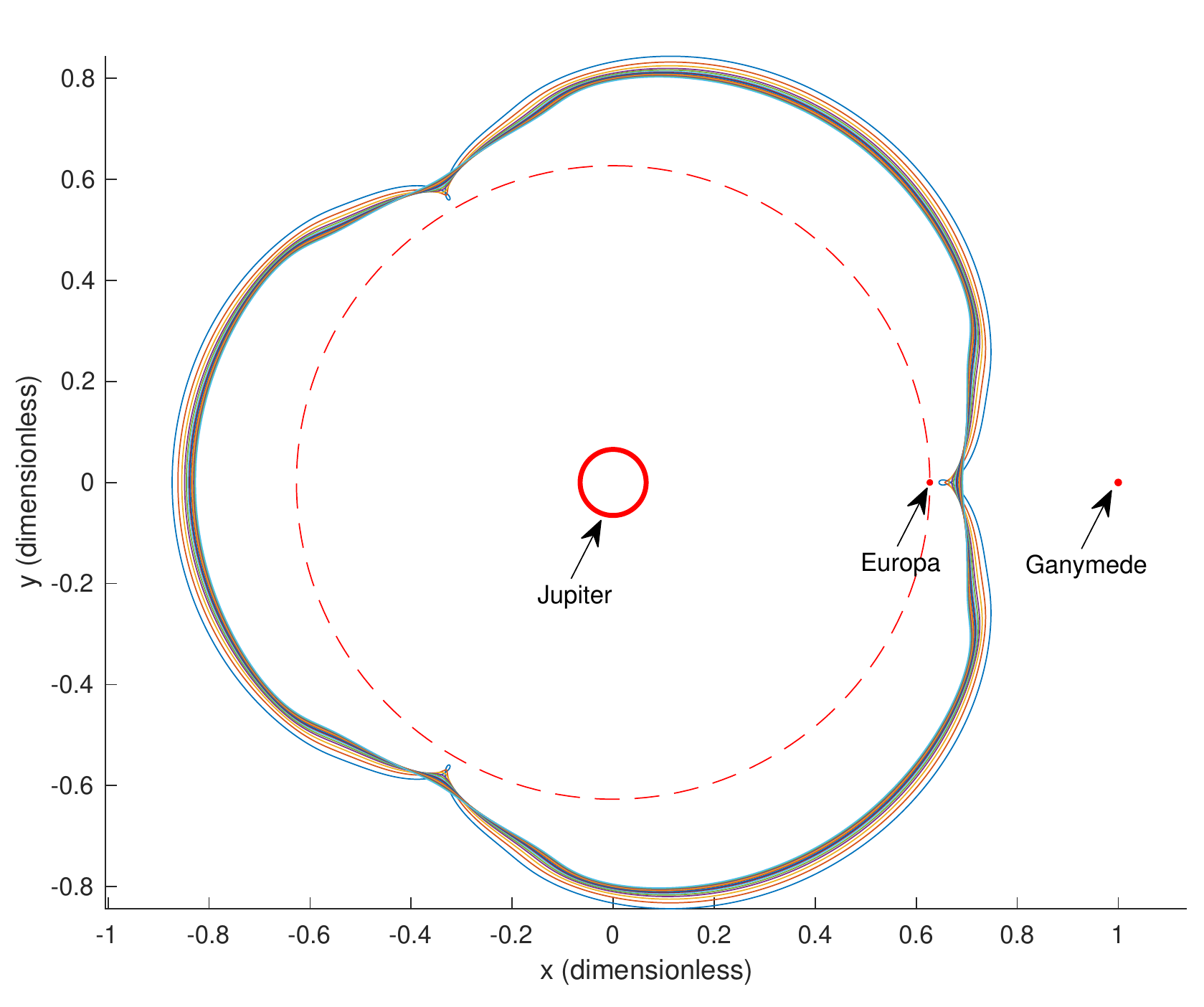}
\caption{ \label{fig:34_1e-5omega} Continuation by $\omega$ of $\bar \mu_{2}=1.0015 \times 10^{-5}$ Jupiter-Europa 3:4 resonant invariant circle from $\omega = 3.111756$ to $\omega = 3.116809$ in Jupiter-Ganymede frame}
\end{centering}
\end{figure}
\begin{figure}
\begin{centering}
\includegraphics[width=0.6\columnwidth]{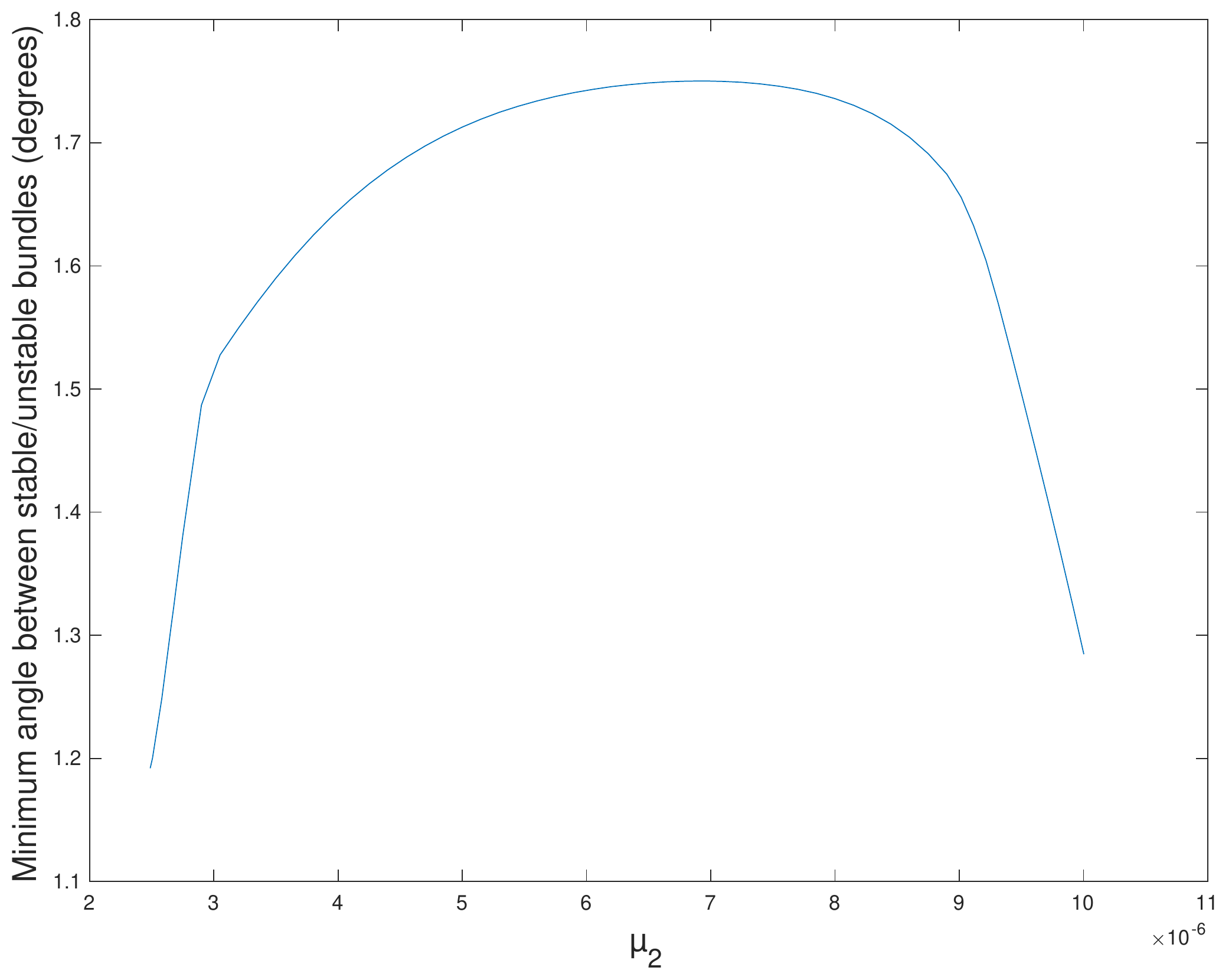}
\caption{ \label{fig:anglevsmu2} Plot of minimum (over $\theta)$ angle between stable and unstable torus directions, versus $\bar \mu_{2}$}
\end{centering}
\end{figure}

 
In Figure \ref{fig:34epsdown}, on the left we show the first part of the continuation by $\bar \mu_{2}$, where we take $\mu_{2}$ from $2.5265 \times 10^{-5}$ (the physical value) to $1.0015 \times 10^{-5}$ for an invariant circle with $\omega = 3.111756$. Since we kept $\omega$ fixed, as $\bar \mu_{2}$ decreases the torus makes closer and closer approaches to Europa; hence, at $\bar \mu_{2}=1.0015 \times 10^{-5}$ we continued the circle by $\omega$ to $\omega = 3.116809$ as shown in Figure \ref{fig:34_1e-5omega}. The resulting torus does not approach Europa as closely, and is less unstable. This circle is then continued to $\bar \mu_{2} = 2.506370 \times 10^{-6}$ (9.9\% of Europa's actual mass) with $\omega$ fixed again; this is shown on the right of Figure \ref{fig:34epsdown}. We once more observe the torus moving closer to Europa as $\bar \mu_{2}$ decreases and $\omega$ is kept fixed. Our quasi-Newton method started diverging when trying to continue the $\omega = 3.116809$ and $\bar \mu_{2} = 2.506370 \times 10^{-6}$ circle further by either $\omega$ or by $\bar \mu_{2}$. We suspect this might be due to a rapid decline, as $\bar \mu_{2}$ decreases to $2.506370 \times 10^{-6}$, in the minimum angle made by the stable and unstable directions to the torus (columns 3 and 4 of $P(\theta)$ from Eq. \eqref{bundleEquations}-\eqref{Lambdaform}). As this angle gets closer to zero, this makes $P$ increasingly close to singular, which negatively affects the convergence of the quasi-Newton method; we plot the minimum angle versus $\bar \mu_{2}$ in Figure \ref{fig:anglevsmu2}. 

\subsection{Structure of Jupiter-Europa 3:4 CCR4BP Torus as $\bar \mu_{2} \rightarrow 0$}

\begin{figure}
\begin{centering}
\includegraphics[width=0.49\columnwidth]{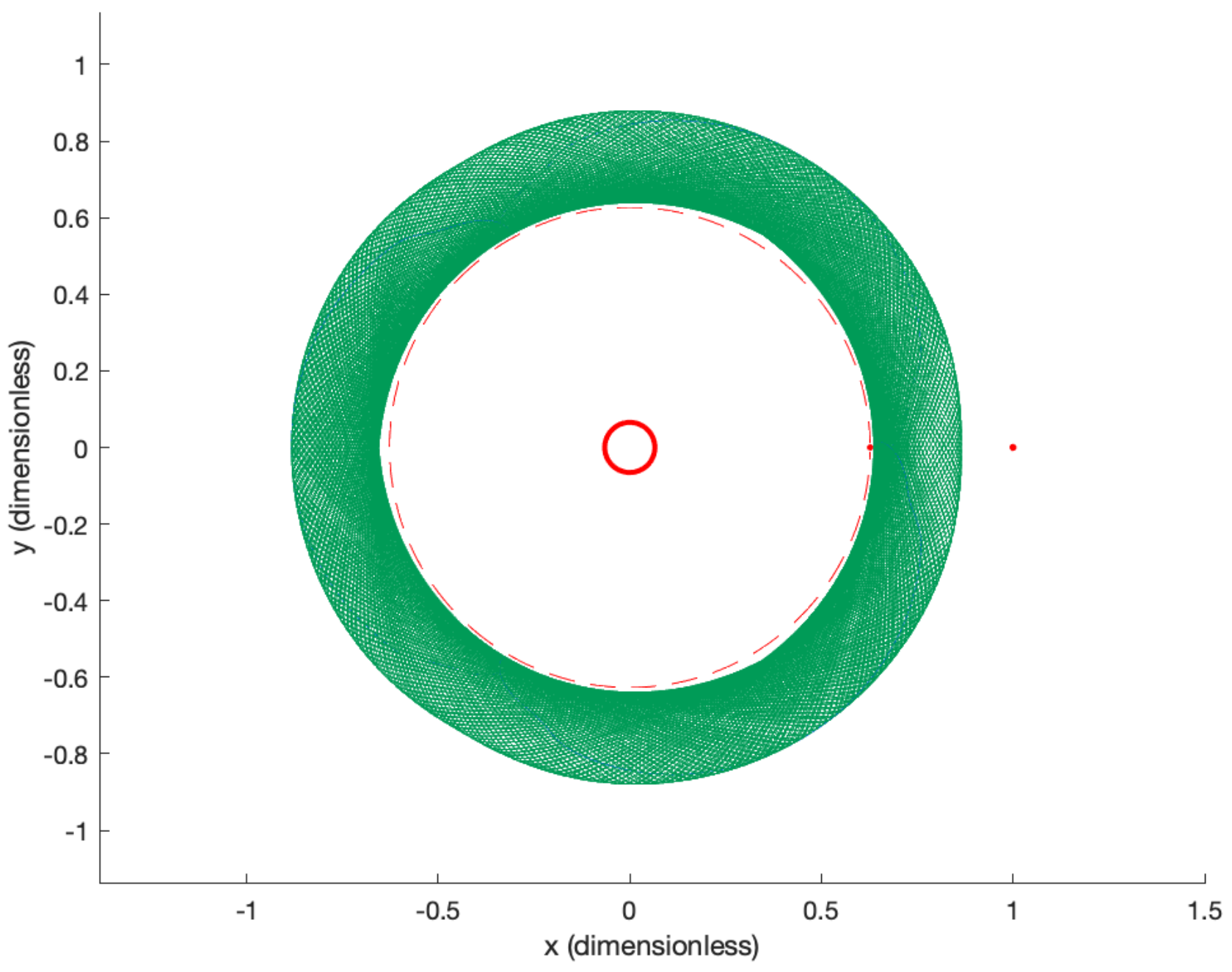}
\includegraphics[width=0.49\columnwidth]{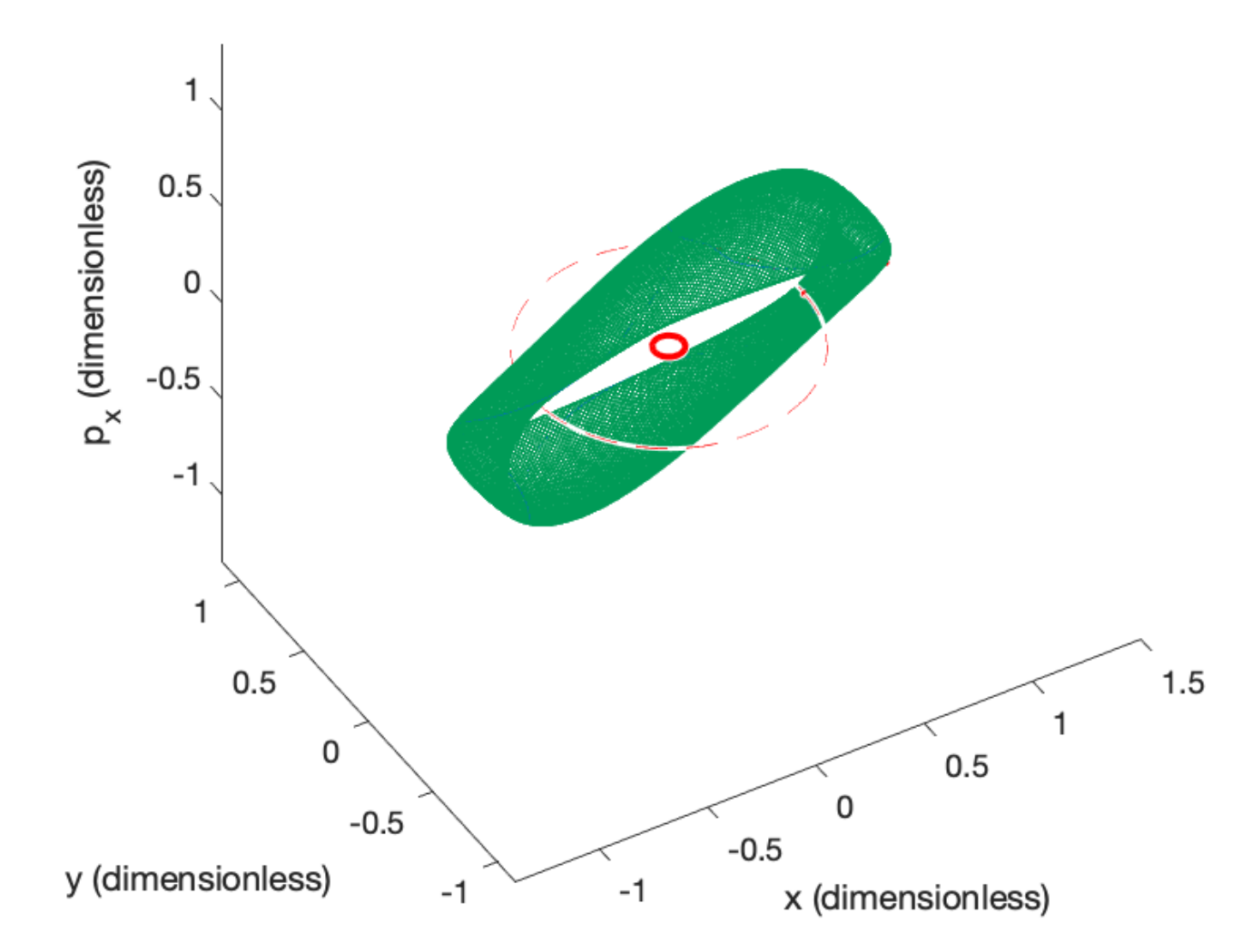}
\caption{ \label{fig:34smalleps} 2D $\bar \mu_{2} = 2.506370 \times 10^{-6}$ CCR4BP Flow-Invariant Torus for Jupiter-Europa 3:4 Resonance in Jupiter-Ganymede Rotating Frame}
\end{centering}
\end{figure}

To characterize the type of Jupiter-Ganymede PCRTBP dynamical structure we may be approaching as $\bar \mu_{2}$ decreases to zero, we integrated trajectories in the Jupiter-Ganymede frame from the $\omega = 3.116809$ and $\bar \mu_{2} = 2.506370 \times 10^{-6}$ invariant circle. As mentioned earlier, this gives us the full flow-invariant 2D torus corresponding to the stroboscopic map-invariant circle. We plotted this 2D torus in both $(x,y)$ space as well as in $(x,y,p_{x})$ space; the results are given in Figure \ref{fig:34smalleps}. Clearly, the resulting structure is still similar to that of the Jupiter-Europa 3:4 resonant CCR4BP tori computed for full $\bar \mu_{2}$, such as the one from Figure \ref{fig:34flowtori}. Our initial hypothesis before carrying out the computations was that by continuing $\bar \mu_{2}$ to zero, we might approach a Jupiter-Ganymede PCRTBP stable 3:2 resonant periodic orbit. However, we see no indications of this. Instead, given the results, what we now believe is occurring is that we are approaching a Jupiter-Ganymede PCRTBP stable 2D non-resonant KAM torus as $\bar \mu_{2} \rightarrow 0$. 

To see this, we converted the torus states shown in Figure \ref{fig:34smalleps} to Jupiter-Ganymede synodic Delaunay coordinates. For full details of this coordinate transformation see Celletti\cite{celletti}; we summarize the basics here. For any state in the Jupiter-Ganymede rotating cartesian coordinate frame, we can compute the instantaneous semi-major axis $a$, eccentricity $e$, longitude of periapse $g$ \emph{relative to the Jupiter-Ganymede rotating $x$-axis}, and the mean anomaly $\ell$. The Delaunay coordinates are comprised of $L = \sqrt{\mathcal{G}m_{1}a}$, $G = L\sqrt{1-e^{2}}$, $\ell$, and $g$. The most important thing about synodic Delaunay coordinates is that they are action-angle variables for the rotating frame Kepler problem (the PCRTBP with $\mu=0$). Hence, in the Kepler problem, the Hamiltonian becomes a function of only $L$ and $G$; the actions $L$ and $G$ remain constant along trajectories, while $\ell$ and $g$ advance at constant frequencies which are a function of $L$ only ($\dot g = -1$, as that is the rotation rate of the coordinate frame). For any $(L,G)$ such that $\dot \ell$ and $\dot g$ are non-resonant, this results in an invariant 2D torus with $(\ell,g)$ values in $\mathbb{T}^{2}$ is densely filled over time by any trajectory starting on the torus.   

When the rotating-frame Kepler problem is perturbed by $\mu>0$ in the CRTBP, the problem is no longer globally integrable; that is, there are no global action-angle coordinates. Nevertheless, many of the 2D stable   non-resonant tori described just earlier persist with just small changes. These tori are called KAM tori, and are discussed further by Celletti\cite{celletti}; indeed, Celletti and Chierchia\cite{celletti2007kam} rigorously proved the existence of such tori inside an energy surface in the Sun-Jupiter CRTBP. Since KAM tori do not change much compared to the tori of the same frequencies in the Kepler problem, the values of the actions $L$ and $G$ do not change significantly on the torus either. However, the full range of possible $(\ell,g) \in \mathbb{T}^{2}$ should still be densely filled in such a KAM torus. There can exist other tori in the perturbed system which do not visit the whole range of all possible $(\ell,g)$ values; these are referred to as secondary tori, and are topologically different from KAM tori. 

\begin{figure}
\begin{centering}
\includegraphics[width=0.49\columnwidth]{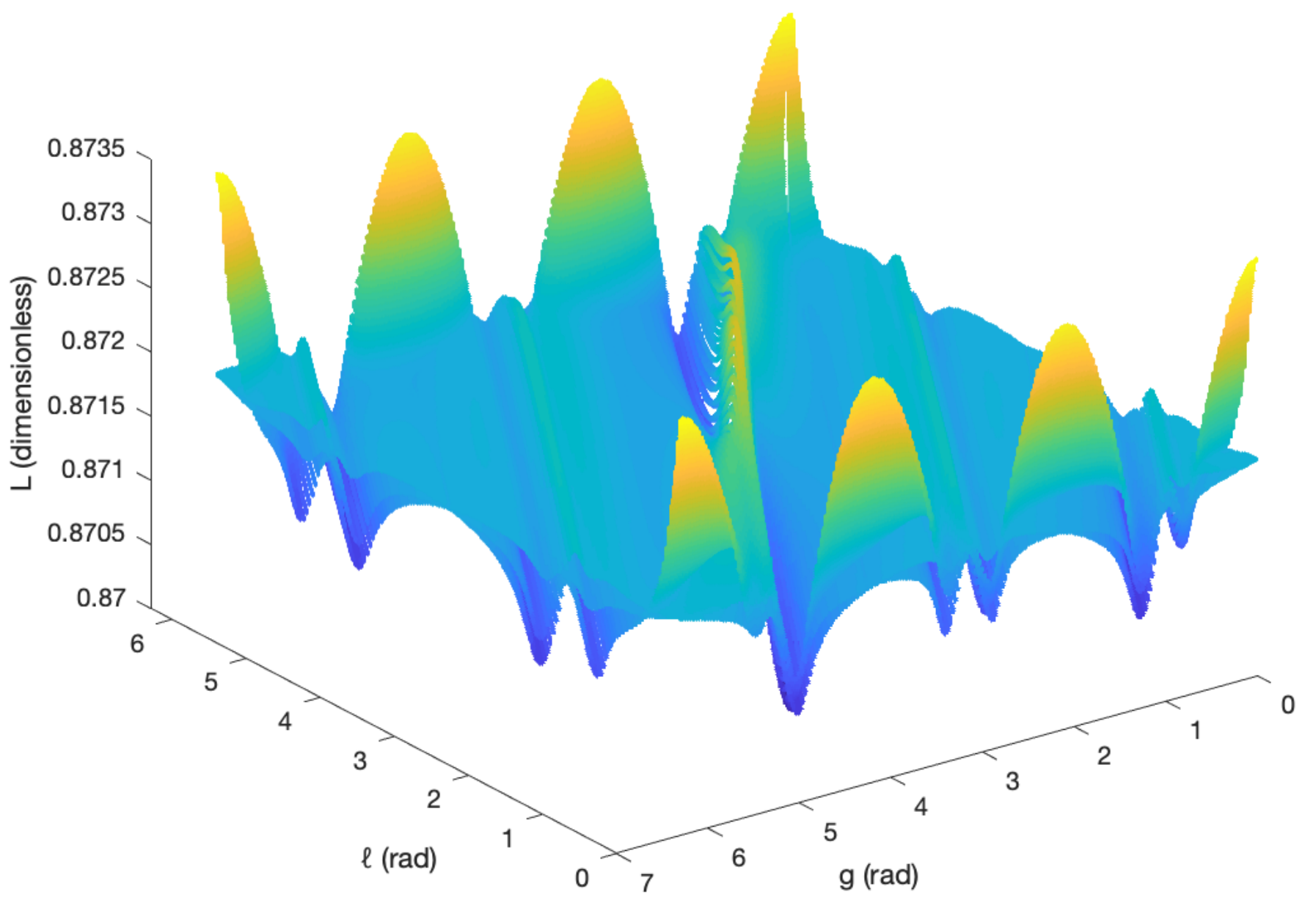}
\includegraphics[width=0.49\columnwidth]{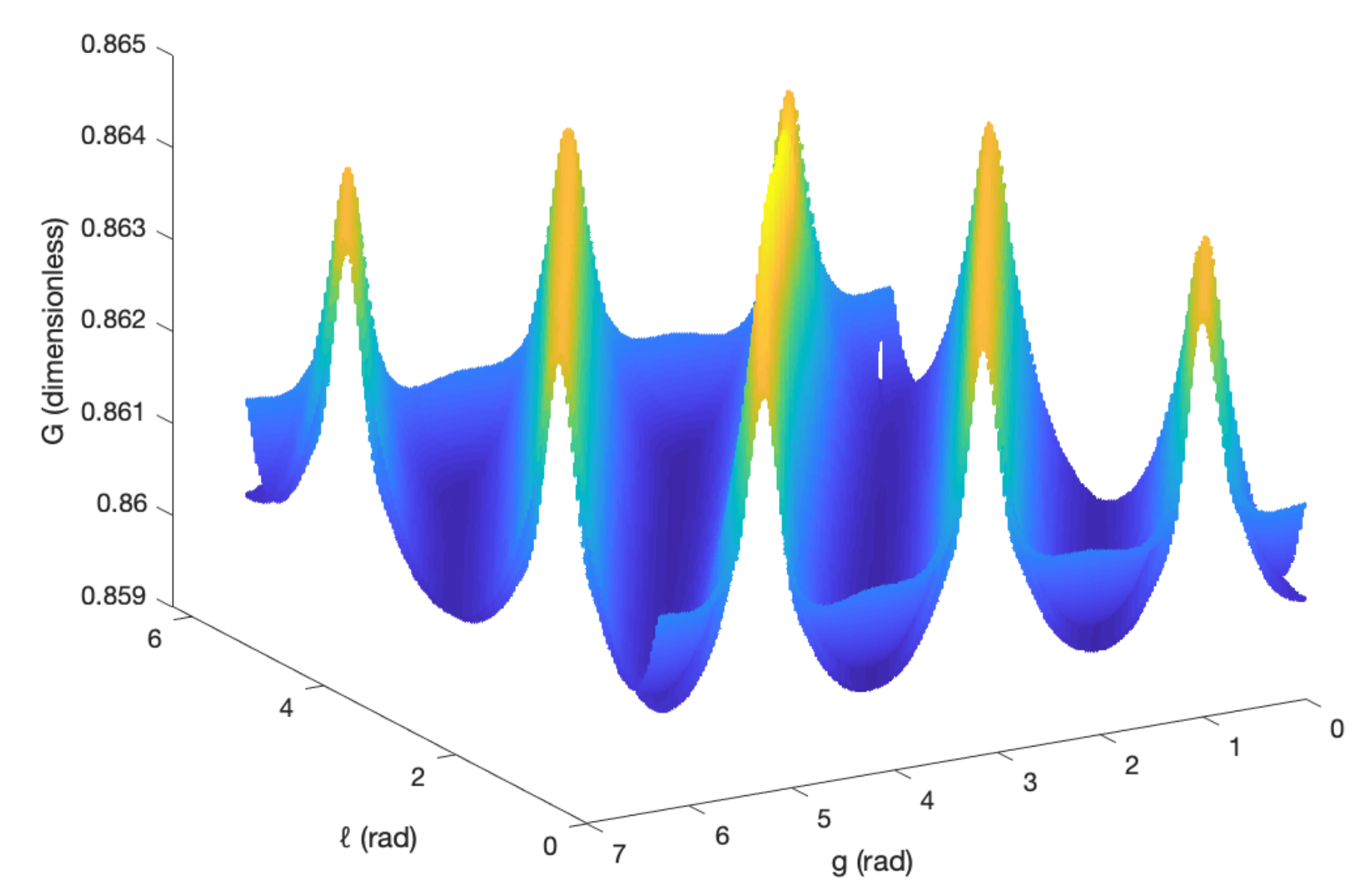}
\caption{ \label{fig:delaunayPlot} 2D $\bar \mu_{2} = 2.506370 \times 10^{-6}$ CCR4BP Flow-Invariant Torus for Jupiter-Europa 3:4 Resonance in Jupiter-Ganymede Synodic Delaunay Coordinates}
\end{centering}
\end{figure}

Converting the torus states shown in Figure \ref{fig:34smalleps} to $(L,G,\ell,g)$ coordinates and plotting the resulting torus in these Delaunay coordinates, it is clear that the torus topology is similar to that expected for a KAM torus. The Delaunay coordinate plots in $(\ell,g,L)$ space and $(\ell,g,G)$ space are shown on the left and right respectively of Figure \ref{fig:delaunayPlot}. We indeed observe that the whole range of $(\ell,g)$ values in $\mathbb{T}^{2}$ is densely filled, while the actions $(L,G)$ do not significantly vary on the torus. Hence, we propose that the limiting dynamical structure of the Jupiter-Europa 3:4 resonant CCR4BP torus in the Jupiter-Ganymede PCRTBP is a non-resonant KAM torus.

We believe that the reason for this is the slight incommensurability of the periods of Europa and Ganymede around Jupiter; Europa and Ganymede are in an approximate, not exact, 2:1 resonance with each other. Hence, writing $\Omega_2$ for the rotational frequency of Europa around Jupiter in the Jupiter-Ganymede frame time units, the frequencies involved in a Jupiter-Europa 3:4 resonant orbit are $3\Omega_2/4$ (mean motion), $-1$ (due to $\dot g$ caused by the frame rotating with Ganymede), and $\Omega_{2}-1$ (due to Europa). This creates a single resonance as $4(3\Omega_2/4)+3(-1)-3(\Omega_{2}-1)=0$; $\Omega_2$ is not rational so there are no other resonances. If Europa's mass is set to zero, the $\Omega_{2}-1$ frequency is decoupled from the other two, so the resonance is broken. The only coupled frequencies left are the $3\Omega_2/4$ and $-1$, which are incommensurate, leading to a KAM torus. 

\section{Conclusions}

Analyzing models incorporating more than two large bodies, such as the CCR4BP, can help to more accurately link CRTBP models corresponding to different pairs of bodies as compared to simpler patched approximations. An understanding of the resonances present in the CCR4BP can provide insight into potentially important structures for multi-moon tour design. In this study, we successfully computed the quasi-periodic analogues of Jupiter-Europa 3:4 and Jupiter-Ganymede 3:2 and 7:5 resonant periodic orbits in the Jupiter-Europa-Ganymede planar CCR4BP. A predictor based on the Poincar\'e-Lindstedt method was coupled with an efficient quasi-Newton method corrector to compute the tori and their stable and unstable directions and multipliers. Significant differences were observed in the shapes of the full 2D flow-invariant tori depending on whether a Jupiter-Europa or Jupiter-Ganymede rotating frame was used. After computing Jupiter-Europa 3:4 tori in the CCR4BP with physical Jupiter, Europa, and Ganymede masses, we found that continuing the torus with Europa mass decreasing towards zero results in the unstable Jupiter-Europa resonant torus likely approaching a stable Jupiter-Ganymede PCRTBP nonresonant KAM torus. 

Through the computations carried out in this paper, we were also able to verify the performance of the quasi-Newton method developed in our previous work\cite{kumar2021rapid}, implementing it for a different dynamical model then the planar elliptic RTBP model used in that paper for numerical demonstrations. The computational results were very positive in this application as well. The tori computed here should prove useful for future work on calculations of low-energy trajectories between them. Indeed, the quasi-Newton method we used for the computations also gave us the linear manifold approximations as columns 3 and 4 of the $P(\theta)$ matrix. We expect to study these manifolds further in the very near future, with a view to searching for heteroclinic connections in the CCR4BP linking Jupiter-Ganymede resonances to Jupiter-Europa resonant orbits. 

\section{Acknowledgments}

This work was supported by a NASA Space Technology Research Fellowship. This research was carried out in part at the Jet Propulsion Laboratory, California Institute of Technology, under a contract with the National Aeronautics and Space Administration (80NM0018D0004). R.L was supported in part by NSF grant DMS 1800241.

\bibliographystyle{AAS_publication}   
\bibliography{references}   

\end{document}